\documentclass[a4paper]{amsart}
\usepackage[final]{graphicx}
\usepackage{a4wide}
\usepackage[english]{babel}
\usepackage[T1]{fontenc}
\usepackage{subcaption}
\usepackage{dsfont}

\usepackage{float}
\usepackage{subfig}
\usepackage{mathrsfs}
\usepackage{amsmath}
\usepackage{amssymb}
\usepackage{verbatim}
\usepackage{amsthm}
\usepackage[usenames]{color}
\usepackage{fancybox}
\usepackage{enumitem}
\usepackage{eurosym}
\usepackage[authoryear]{natbib}
\usepackage{rotating}
\def\N{\mathbb{N}}
\def\Pro{\mathbb{P}}
\def\E{\mathbb{E}}
\def\R{\mathbb{R}}
\def\1{\large{1}}
\def\Acal{\mathcal{A}}
\def\Lcal{\mathcal{L}}
\def\Fcal{\mathcal{F}}
\def\Hcal{\mathcal{H}}
\def\Rcal{\mathcal{R}}
\newcommand{\Fourier}[1]{\widetilde{#1}}%
\newcommand{\paren}[1]{\left( #1 \right)} %
\newcommand{\accolade}[1]{\left\{ #1 \right\}}%

\newcommand{\FigErrPlots}[3]{\includegraphics[width=.4\linewidth]{#1-#3-#2}}
\newcommand{\FigPlotsLS}[3]{\includegraphics[width=.23\linewidth]{#1-rho-#2#3ls}}
\newcommand{\FigPlotsSURF}[3]{\includegraphics[width=.36\linewidth]{#1-rho-#2#3surf}}

\theoremstyle{plain}
\newtheorem{Theorem}{Theorem}
\newtheorem{Lemma}{Lemma}
\newtheorem{Proposition}{Proposition}
\newtheorem{Corollary}{Corollary}

\graphicspath{{./images/}}
\usepackage{url}

\title[Minimax and adaptive estimation of dthe Wigner function in QHT]{Minimax and adaptive estimation of the Wigner function in quantum homodyne tomography with noisy data}

\begin{document}
\maketitle

\begin{center} {Karim Lounici$^{1}$ ,
Katia Meziani$^{2}$
and 
Gabriel Peyr\'e$^{3}$}\end{center}
\begin{center} CEREMADE, UMR CNRS 7534, Universit\'e Paris Dauphine $^{2,3}$,\\School of Mathematics, Georgia Inst. of Technology$^{1}$, \\  CREST-ENSAE$^{2}$\end{center}
{\textit{klounici@math.gatech.edu}, \textit{meziani@ceremade.dauphine.fr}, \textit{gabriel.peyre@ceremade.dauphine.fr} }

\begin{abstract} ~In quantum optics, the quantum state of a light beam  is represented through the Wigner function, a density on $\R^2$ which may take negative values but must respect intrinsic positivity constraints imposed by quantum physics. In the framework of noisy quantum homodyne tomography with efficiency parameter $1/2 < \eta \leq 1$, we study the theoretical performance of a kernel estimator of the Wigner function. We prove that it is minimax efficient, up to a logarithmic factor in the sample size, for the $\mathbb L_\infty$-risk over a class of infinitely differentiable. We compute also the lower bound for the $\mathbb L_2$-risk. We construct adaptive estimator, i.e. which does not depend on the smoothness parameters, and prove that it attains the minimax rates for the corresponding smoothness class
functions. Finite sample behaviour of our adaptive procedure are explored through numerical experiments.
    \end{abstract}

\begin{center}
\end{center}



\smallskip

\noindent  This paper deals with a \textit{severely ill-posed inverse problem} which comes from quantum optics. Quantum optics is a branch of quantum mechanics which studies  physical systems at the atomic and subatomic scales. Unlike classical mechanics, the result of a physical measurement is generally random. Quantum mechanics does not predict a deterministic course of events, but rather the probabilities of various alternative possible events. It provides predictions on the outcome measures, therefore explore measurements  involve non-trivial statistical methods and inference on the result of a measurement should to be done  on identically prepared quantum system.  \\
\noindent To understand our statistical model, we start in Section~\ref{physbackground} with a short introduction to the needed quantum notions. Section~\ref{Statistical.model} introduces the statistical model by making the link with quantum theory. Interested reader can get further acquaintance with quantum concepts through the textbooks  or the review articles of \cite{Helstrom,Holevo,Barndorff-Nielsen&Gill&Jupp} and \cite{Leonhardt}.

\section{Physical background}
\label{physbackground}
 In quantum mechanics, the measurable properties (ex: spin, energy, position, ...) of a quantum system are called "\textit{observables}". The probability of obtaining each of the possible outcomes when measuring an observable is encoded in the \textit{quantum state} of the considered physical system. 
 
\subsection{Quantum state and observable}

\noindent  The mathematical description of the quantum state of a system is
given in form of a density operator  $\rho$ on a complex Hilbert space $\Hcal$ (called the space of states) satisfying the three following conditions:
\begin{enumerate}
     \item Self adjoint: $\rho=\rho^{*}$, where $\rho^{*}$ is the adjoint of $\rho$.
     \item Positive: $\rho\geq 0$, or equivalently $\langle\psi,\rho\psi\rangle\geq0$ for all
      $\psi\in \mathcal  {H}$.
     \item Trace one: $\mathrm{Tr}(\rho)=1$.
\end{enumerate}
Notice that $\mathcal D(\mathcal H) $ the set of density operator $\rho$ on $\mathcal H$ is a convex set.  The extreme points of the convex set $\mathcal D(\mathcal H) $ are called \textit{pures states} and all others states are called \textit{mixed states}.\\

\begin{figure}[!lh]
\centering
\caption{The density matrix $\rho$ of a coherent-3 state.}
\includegraphics[scale=0.5]{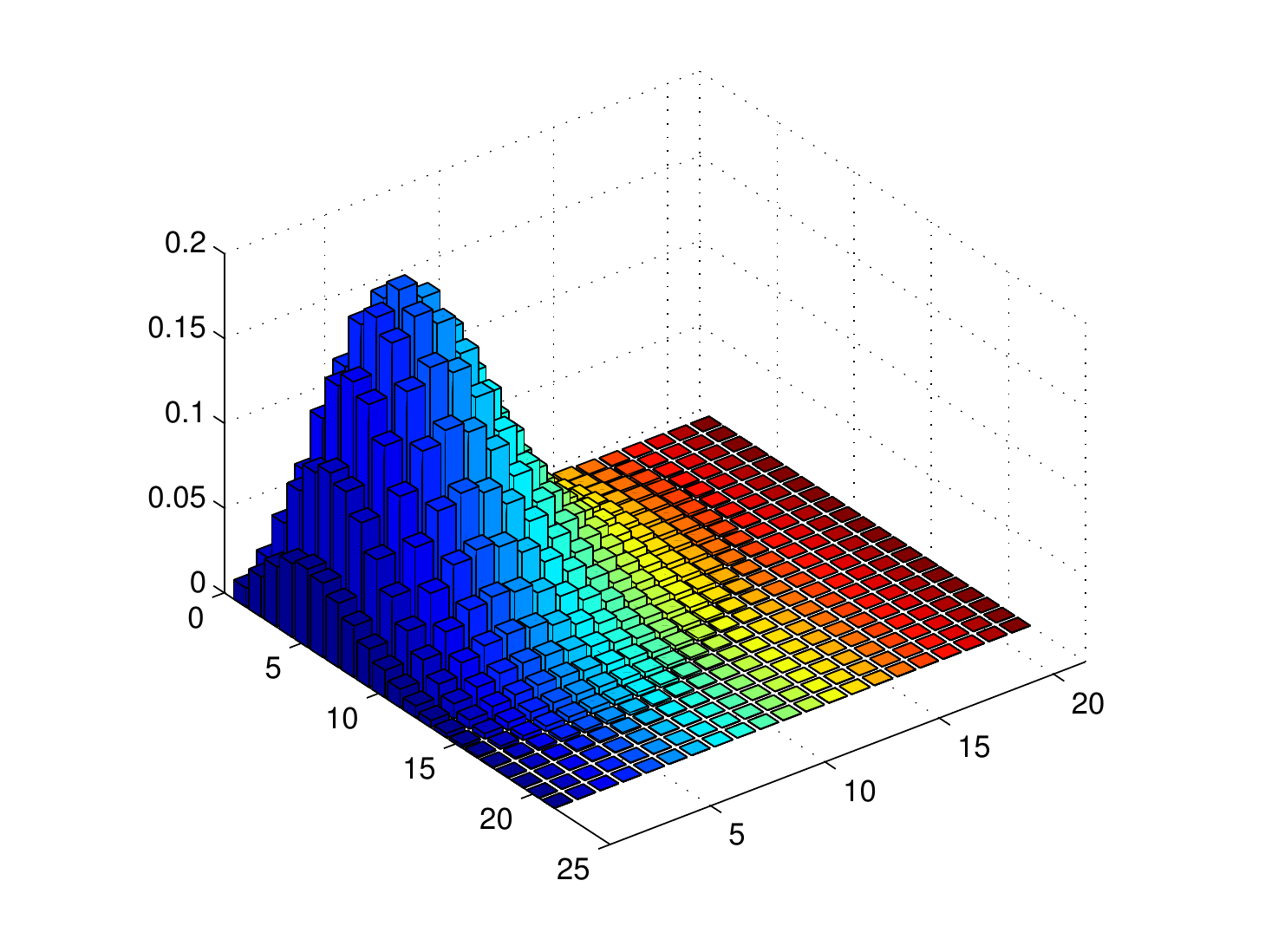}
\label{coherent3}
\end{figure}

\noindent   In this paper, the quantum system we are interested in is a monochromatic light in a cavity. In this setting of quantum optics, the
space of states $\Hcal$ we are dealing with is the space of square integrable
complex valued functions on the real line.  A particular orthonormal basis comes with this Hilbert space is the Fock basis $\accolade{\psi_j}_{j\in\N}$:
\begin{equation}
          \label{eq:fock}
          \psi_j(x):= \frac{1}{\sqrt{\sqrt{\pi}2^jj!}}H_j(x)e^{-x^2/2},
\end{equation}
where $H_j(x) := (-1)^j e^{x^2} \frac{d^j}{dx^j} e^{-x^2}$ denote the $j$-th Hermite polynomial. In this basis, a quantum state is described  by an infinite density matrix $\rho=[\rho_{j,k}]_{j,k\in\N}$ whose entries are equal to
$$
\rho_{j,k}=\langle \psi_j,\rho\psi_k\rangle,
$$
with $\langle \cdot,\cdot\rangle$ the inner product. The quantum states which can be created at this moment in laboratory are matrices whose entries are decreasing exponentially to 0, \textit{i.e.}, belong to the natural class $\Rcal(C,B,r)$ defined bellow, with $r=2$. Let us define for $C\geq 1$, $B>0$ and
$0 < r \leq 2$, the class $\Rcal(C,B,r)$ is defined as follow
\begin{equation}
        \label{eq.classcoeff}
         \Rcal(C,B,r) := \{\rho {\rm \ quantum \ state} :
         |\rho_{m,n} |\leq C \exp(-B (m+n)^{r/2})\}.
\end{equation}
\noindent   An example of density matrix of a pure state whose entries are real is given in Figure~\ref{coherent3}.

\smallskip
\noindent  In order to describe mathematically a measurement performed on an observable of a quantum system prepared in state $\rho$, we give the mathematical description of an observable. An observable $\mathbf{X}$ is a self adjoint operator on the same space
of states $\Hcal$ and
$$
\mathbf{X}=\sum_{a}^{dim\Hcal}x_a\mathbf{P}_a,
$$
where the eigenvalues $\{x_a\}_a$ of the observable $\mathbf{X}$ are real and  $\mathbf{P}_a$  is the projection onto the one dimensional space generated by
the eigenvector of $\mathbf{X}$ corresponding to the eigenvalue $x_a$. \\
\noindent  As a quantum state $\rho$ encompasses all the probabilities of the observables of the considered quantum system, when performing a measurement of the observable $\mathbf{X}$ of a quantum state $\rho$,  the result is a random variable $X$ with values in the set of the eigenvalues of the observable $\mathbf{X}$. For a quantum system prepared in state $\rho$,  $X$ has the following probability distribution and expectation function
\begin{center}
$\Pro_\rho(X=x_a)=\mathrm{Tr}(\mathbf{P}_a\rho)\quad$ and
$\quad\E_\rho(X)=\mathrm{Tr}(\mathbf{X}\rho)$.
\end{center}
\noindent  Note that the conditions defining the density matrix $\rho$ insure that $\Pro_\rho$ is a probability distribution. In particular, the characteristic function is given by
$$
\E_\rho(e^{itX})=\mathrm{Tr}(\rho e^{it\mathbf{X}}).
$$

\subsection{Quantum homodyne tomography and Wigner function}
\label{QHTWfunction}

\noindent In quantum optics, a monochromatic light in a cavity is described by a quantum harmonic oscillator. In this setting, the observables of interest are usually $\textbf{Q}$ and $\textbf{P}$ (resp. the electric and magnetic fields). But according to Heisenberg's uncertainty principle,  $\textbf{Q}$ and $\textbf{P}$ are non-commuting observables,  they may not be simultaneously measurable. Therefore, by performing measurements on $(\textbf{Q},\textbf{P})$, we cannot get a probability density of the result $(Q,P)$. However, for all phase $\phi\in[0,\pi]$ we can measure the quadrature observables 
 $$
 \mathbf{X}_\phi :=\mathbf{Q}\cos \phi +\mathbf{P}\sin \phi.
 $$
\noindent  Each of these quadratures could be measured on a laser beam by a technique put in practice for the first time by \textit{Smithey} and called
\textbf{Quantum Homodyne Tomography} (QHT). The theoretical foundation of quantum homodyne tomography was outlined by \cite{Vogel&Risken}.

\begin{figure}[!lh]
\centering
\caption{QHT measurement scheme.}
\includegraphics[scale=0.5]{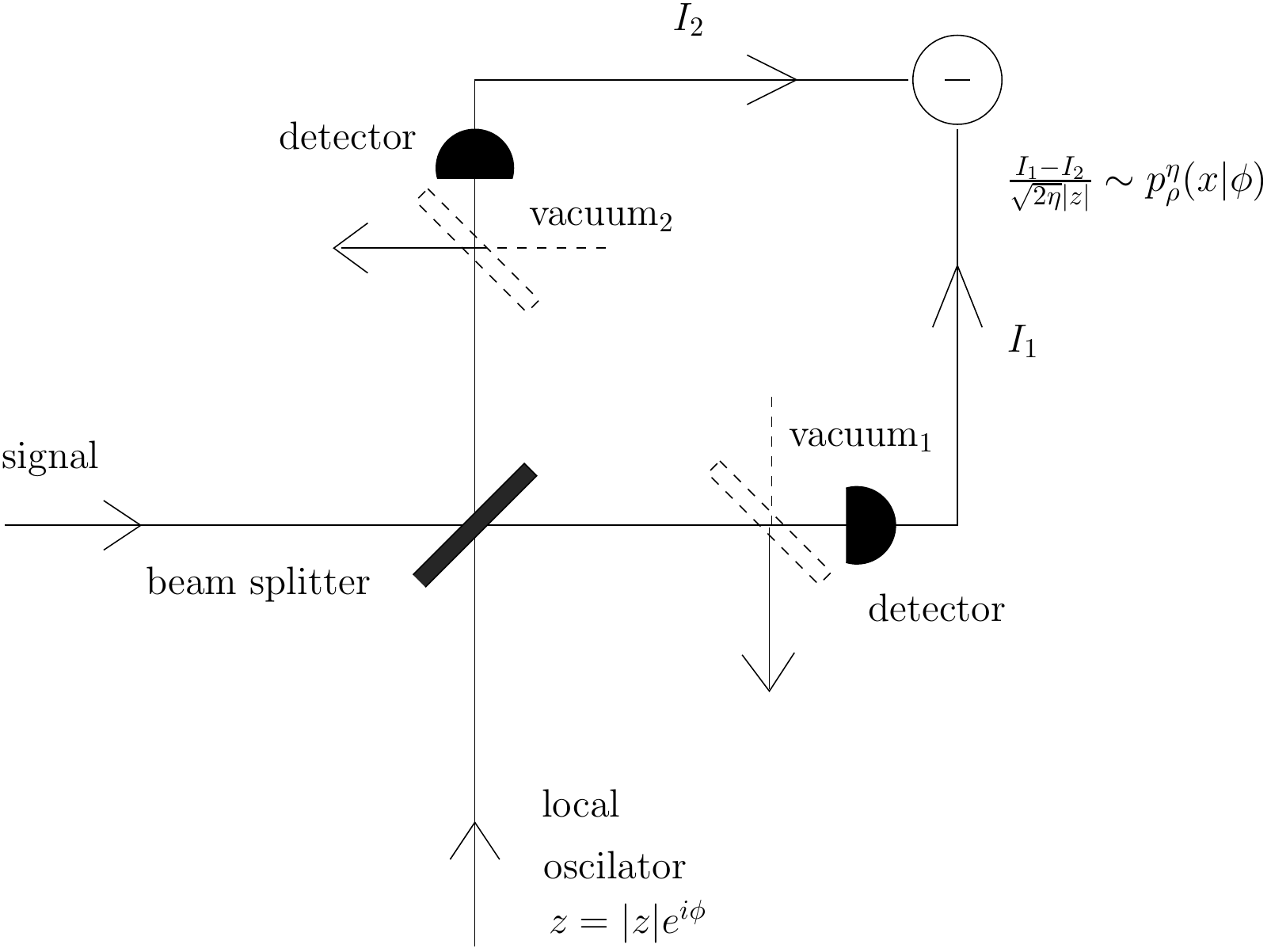}
\label{Pfig:1}
\end{figure}

\smallskip
\noindent The experimental set-up, described in Figure~\ref{Pfig:1}, consists of mixing the signal field  with a local oscillator field (LO) of high intensity $\left|z\right|>>1$. The phase $\Phi$ of the LO is choosen s.t. $\Phi\sim\mathcal{U}[0,\pi]$. The resulting beam is split by a 50-50 beam splitter, and the photodetectors count the photons in the two output beams by giving integrated currents $I_1$ and $I_2$ proportional to the number of photons. The result of the measurement is
produced by taking the difference of the two currents and rescaling it by the intensity $|z|$.   In the case of noiseless measurement and for a phase $\Phi=\phi$, the result $X_\phi=\frac{I_2-I_1}{|z|}$ has density $p_\rho(\cdot|\phi)$ corresponding to measuring $\mathbf{X}_\phi$.\\

\noindent  In others words, when performing a QHT measurement of the observable $\mathbf{X}_\phi$ of the quantum state $\rho$,  the result is a random variable $X_\phi$ whose density conditionally to $\Phi=\phi$ is denoted by $p_\rho(\cdot|\phi)$.  It's characteristic function is given by
$$
\E_\rho(e^{itX_\phi})=\mathrm{Tr}(\rho e^{it\mathbf{X}_\phi})=\mathrm{Tr}(\rho e^{it(\mathbf{Q}\cos \phi +\mathbf{P}\sin \phi)})=\mathcal F_1[ p_\rho(\cdot|\phi) ](t),
$$
where $\mathcal F_1[ p_\rho(\cdot|\phi) ](t)=\int e^{itx}p_\rho(\cdot|\phi)dx$ denotes the Fourier transform with respect to the first variable. Moreover if $\Phi$ is chosen uniformly on $[0,\pi]$, the joint density probability of $(X_\phi,\Phi)$ with respect to the Lebesgue measure on $\R\times[0,\pi]$ is
$$
p_\rho(x,\phi)=\frac 1 \pi p_\rho(x|\phi)\1_{[0,\pi]}(\phi).
$$

\noindent  An equivalent representation for a quantum state $\rho$ is the function $W_\rho:\R^2\rightarrow\R$ called the Wigner function, introduced for the first time by \cite{Wigner}.  The Wigner function may be obtained from the momentum representation
\begin{eqnarray}
          \hspace{-1cm}
         \label{def.fourierWigner}        
          \widetilde{W}_{\rho}(u,v)\hspace{-0.15cm}&:=&\hspace{-0.15cm}\Fcal_2[W_\rho]
          (u,v)=\mathrm{Tr}\left(\rho e^{i(u\textbf{Q}+v  \textbf{P})}\right),
 \end{eqnarray}
where $\Fcal_2$ is its Fourier transform with respect to both variables. By applying a change of variables $(u,v)$ into $(t\cos \phi,t\sin\phi)$, we get 
\begin{eqnarray}
            \hspace{-1cm}
            \label{def.fourierRadon}        
             \widetilde{W}_{\rho}(t\cos \phi,t\sin\phi)=\mathcal F_1[ p_\rho(\cdot|\phi) ](t)=
             \mathrm{Tr}(\rho e^{it\mathbf{X}_\phi}).
 \end{eqnarray}
\noindent  The origin of the appellation quantum homodyne tomography comes from the fact that the procedure described above is similar to positron emission tomography (PET), where the density of the  observations is the Radon transform of the underlying distribution
\begin{equation}
         \label{eq.Radon.transform}
          p_\rho(x|\phi)=\Rcal[W_{\rho}](x,\phi)=\int W_{\rho}(x\cos\phi + t\sin\phi, \, x \sin
          \phi - t\cos\phi )dt,
\end{equation}
\noindent  where $\Rcal[W_{\rho}]$ denotes the Radon transform of $W_\rho$. The main difference with PET is that the role of the unknown distribution is played by the Wigner function which can be negative. \\
\noindent  The physicists consider the Wigner function as a quasi-probability density of $(Q,P)$ if one can measure simultaneously $(\textbf{Q},\textbf{P})$. Nevertheless, the Wigner function does not satisfy all the properties of a conventional probability density but satisfies boundedness properties unavailable for classical densities. For instance, the Wigner function can and normally does go negative for states which have no classical model. The Wigner function is such that
\begin{eqnarray}
       \hspace{-1cm}
        \label{int1}
        W_\rho:\R^2\rightarrow\R,\quad \iint W_\rho(q,p)dqdp=1.
 \end{eqnarray}
\noindent Therefore,  the negative part of the Wigner function makes the interpretation in term of density of probability in space phases less intuitive.  However, the Radon transform  of  the Wigner function is always a probability density. Indeed, conditionally to $\Phi=\phi$ and by applying the change of variables $(q,p)$ into $(x\cos\phi + t\sin\phi, \, x \sin\phi - t\cos\phi )$, it comes
\begin{eqnarray*}
\hspace{-1cm}
               1\hspace{-0.15cm}&=&\hspace{-0.15cm}\iint W_\rho(q,p)dqdp\\
               \hspace{-0.15cm}&=&\hspace{-0.15cm}\iint W_{\rho}(x\cos\phi + t\sin\phi, \, x
                \sin\phi - t\cos\phi )dtdx\\
                \hspace{-0.15cm}&=&\hspace{-0.15cm}\int \Rcal[W_{\rho}](x,\phi)dx=
                \int p_\rho(x|\phi)dx.
 \end{eqnarray*}

\begin{figure}[!h]
  \centering
\includegraphics[scale=0.5]{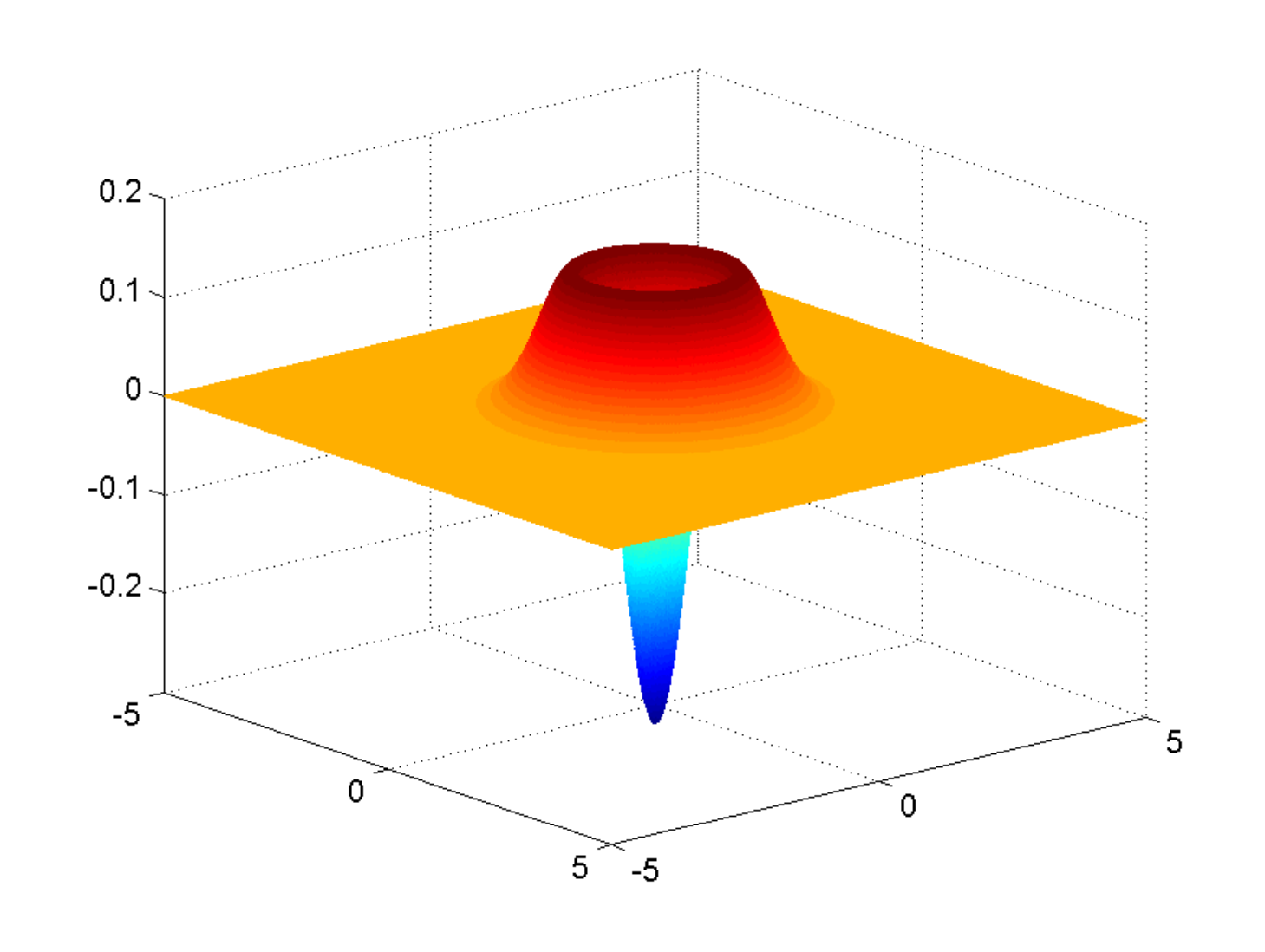}
\includegraphics[scale=0.5]{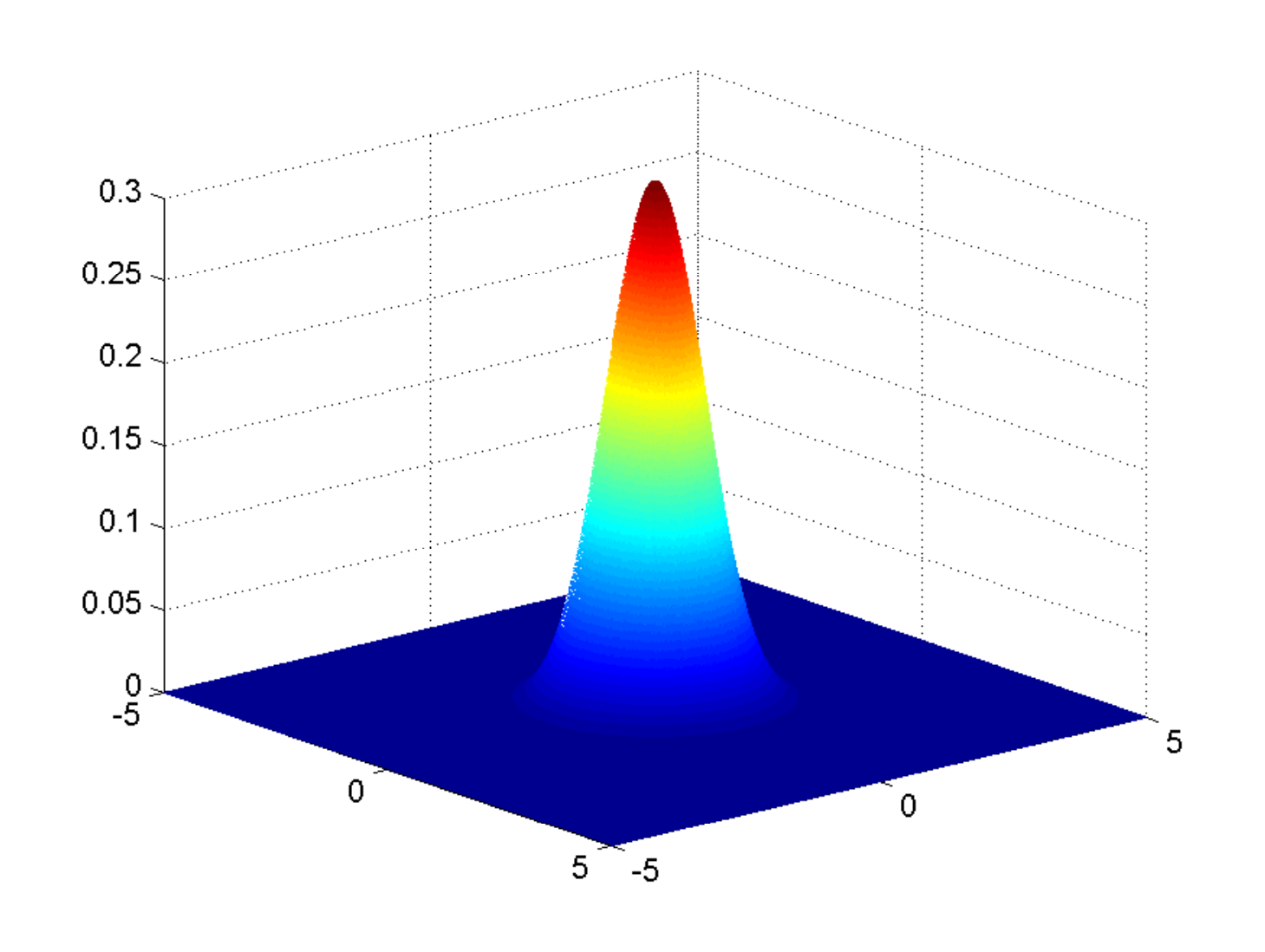}
\caption{The Wigner function of the single photon state (left) and the Wigner function of the vacuum state (right).}
\label{wigner1photon}
\end{figure}

\noindent  Note that  the existence of negative values in the function of Wigner  can be precisely taken like criterion to discriminate nonclassical states of the field. Figure~\ref{wigner1photon} is the representation of the Wigner function of the vacuum state and the nonclassical one-photon state.

\noindent In the Fock basis, we can write $W_{\rho}$ in terms of the density matrix $[\rho_{jk}]$ as follows (see \cite{Leonhardt} for the details).
\begin{equation}
       W_{\rho}(q,p) = \sum_{j,k} \rho_{jk} W_{j,k}(q,p) \nonumber
\end{equation}
where for $j \geq k$,
\begin{eqnarray}
         \hspace{-1cm}
          W_{j,k}(q,p)\hspace{-0.15cm}&=&\hspace{-0.15cm}\frac{(-1)^j}{\pi}
           \paren{\frac{k!}{j!}}^\frac12 \left( \sqrt2(ip-q)\right)^{j-k}e^{-\paren{q^2+p^2}}, 
           L_k^{j-k}\left(2q^2+2p^2\right).
           \label{Weq:WmnLaguerre}
\end{eqnarray}
and $L_{k}^{\alpha}(x)$ the Laguerre polynomial of degree $k$ and order $\alpha$.

\subsection{Pattern functions}
\label{Patternfunction}

\noindent The ideal result of the QHT measurement provide $(X_\phi,\Phi)$ of  joint probability density with respect to the Lebesgue measure on $\R\times[0,\pi]$ equals to
\begin{eqnarray}
         \hspace{-1cm}
          \label{eq:prhophi}
          p_\rho(x,\phi)=\frac 1 \pi p_\rho(x|\phi)\1_{[0,\pi]}(\phi)=\frac 1 \pi \Rcal[W_{\rho}].
          (x,\phi)\1_{[0,\pi]}(\phi)
\end{eqnarray}
\noindent The density $p_\rho(\cdot,\cdot)$  can be written in terms of the entries of the density matrix $\rho$ (see \cite{Leonhardt}) 
\begin{eqnarray}
          \hspace{-1cm}
          \label{de.p.a.rho}
          p_\rho(x,\phi)=\sum_{j,k=0}^\infty \rho_{j,k}\psi_j(x)\psi_k(x)e^{-(j-k)\phi},
\end{eqnarray}
\noindent where $\left\{\psi_j\right\}_{j\in\N}$ is the Fock basis defined in (\ref{eq:fock}). Inversely (see \cite{DAriano.0,Leonhardt} for details), we can write
\begin{eqnarray}
          \hspace{-1cm}
          \label{de.rho.a.p}
          \rho_{j,k}=\int_0^\pi\hspace{-0.25cm}\int p_\rho(x,\phi) f_{j,k}(x)e^{(j-k)\phi}dx
          d\phi,
\end{eqnarray}
\noindent where the functions $f_{j,k}:\R\rightarrow\R$ introduced by \cite{DAriano.3} are called the \textit{"pattern functions"}.  A explicit form of $f_{j,k}(\cdot)$ is given by its Fourier transform by \cite{Richter1}: for all $ j \geq k$
\begin{eqnarray}
          \hspace{-1cm}
           \label{pattern}
            \widetilde{f}_{j,k}(t) \hspace{-0.15cm}&= \widetilde{f}_{k,j}(t)=
            \pi(-i)^{j-k}\sqrt{\frac{2^{k-j} k!}{j!}}|t|
            t^{j-k}e^{-\frac{t^2}{4}}L^{j-k}_k(\frac{t^2}{2}),
\end{eqnarray}
\noindent where $L^{\alpha}_k(x)$ denotes the generalized Laguerre polynomial of degree $k$ and order $\alpha$.  Note that by Writing $t=||w||=||(q,p)||= \sqrt{q^2+p^2} $ in the equation~(\ref{Weq:WmnLaguerre}), we can define
\begin{equation}
         \label{Weq:deflmn}
          l_{j,k}(t) := \left|W_{j,k}(q,p)\right| = \frac{2^{\frac{j-k}2}}{\pi}
          \paren{\frac{k!}{j!}}^\frac12 t^{j-k}e^{-t^2}  \left|L_k^{j-k}(2t^2)\right |.
\end{equation}
Therefore, there exists an useful relation, for all $j\geq k$
\begin{equation}
    \label{pattljk}
\left|\widetilde{f}_{j,k}(t)\right|=\pi^2|t|  l_{j,k}(t/2). 
\end{equation}

\noindent Moreover \cite{ABM} have given the following Lemma which will be useful to prove our main results.

\begin{Lemma}[\cite{ABM}]
\label{ljk}
\noindent\\
\noindent For all $j,k \in \N$ and $J:= j+k+1$, for all $t \geq 0$,
\begin{equation}
                  l_{j,k}(t) \leq \frac1\pi\left\{ \begin{array}{ll}
                                 1 &{\rm \ if \ } 0 \leq t \leq \sqrt{J}, \\
                                 e^{-(t-\sqrt{J})^2} &{\rm \ if \ } t \geq \sqrt{J}.
                                 \end{array} 
                               \right. 
\label{Weq:estimeLaguerre}
\end{equation}
\end{Lemma}

\section{Statistical model}
\label{Statistical.model}

\noindent  In practice, when one performes a QHT measurement (see Figure~\ref{Pfig:1}), a number of photons fails to be detected. These losses may be
quantified by one single coefficient $\eta\in[0,1]$, such that  $\eta=0$ when there is no detection and  $\eta=1$  corresponds to the ideal case (no loss). The quantity $(1-\eta)$ represents the proportion of photons which are not detected due to various losses in the measurement process. The parameter $\eta$ is supposed to be known, as physicists argue, that their machines actually have high detection efficiency, around $0.9=\eta$. In this paper we consider $\eta\in]1/2,1]$.  Moreover, as the detection process is inefficient, an independent gaussian noise interferes additively with the ideal data $X_\phi$.  Note that the gaussian nature of the noise is imposed by the gaussian nature of the vacuum state which interferes additively (see figure~\ref{wigner1photon}).

\noindent To resume,  for $\Phi=\phi$, the effective result of the QHT measurement  is  for a known efficiency $\eta\in]1/2,1]$,  
\begin{equation}
        \label{Y}
         Y=\sqrt{\eta}\,X_\phi+\sqrt{(1-\eta)/2}\,\xi
\end{equation}
\noindent where $\xi$ is a standard Gaussian random variable, independent of the random variable $X_\phi$  having density, with respect to the Lebesgue measure on $\R\times [0,\pi]$, equal to $p_{\rho}(\cdot,\cdot)$ defined in equation~(\ref{eq:prhophi}).  For the sake of simplicity, we re-parametrize (\ref{Y}) as follow
 \begin{equation}
           \label{Z}
           Z:= Y/\sqrt{\eta}=X_\phi+\sqrt{(1-\eta)/(2\eta)}\,\xi:=X_\phi+\sqrt{2\gamma}\,\xi,
\end{equation}
\noindent where $\gamma=(1-\eta)/(4\eta)$ is known and $\gamma\in [0,1/4[$ as $\eta\in]1/2,1]$. Note that $\gamma=0$ corresponds to the ideal case.\\
\noindent Let us denote by $p_{\rho}^{\gamma}(\cdot,\cdot)$ the density of $(Z,\Phi)$ which is the convolution of the density of $X_\phi$
with $N^\gamma(\cdot)$ the density of a centered Gaussian distribution having variance $2\gamma$, that is 
\begin{eqnarray}
          \hspace{-1cm}
          \label{densitbrui}
          p^\gamma_\rho(z,\phi) \hspace{-0.15cm}&=&\hspace{-0.15cm}\left[ \frac 1 \pi 
          \Rcal[W_{\rho}](\cdot,\phi)\1_{[0,\pi]}(\phi)\right] \ast N^\gamma(z)=  p_\rho
          \left(\cdot,\phi\right) \ast N^\gamma(z)\\
           \hspace{-0.15cm}&=&\hspace{-0.15cm} \int p_\rho\left(z-x,\phi \right)
           N^\gamma (x)dx.\nonumber
 \end{eqnarray}

\noindent For $\Phi=\phi$, a useful equation in the Fourier domain, deduced by the previous relation (\ref{densitbrui}) and equation~(\ref{def.fourierRadon}) is
 \begin{equation}
          \label{fourierproun}
          \Fcal_1[p^\gamma_\rho(\cdot,\phi)](t)
          = \Fcal_1[p_\rho(\cdot,\phi)](t)
           \widetilde{N}^\gamma(t)=\Fourier{W}_\rho(t\cos(\phi),t\sin(\phi))
           \widetilde{N}^\gamma(t),
 \end{equation}     
\noindent where $\Fcal_{1} $ denotes the Fourier transform with respect to the first variable and the Fourier transform of $N^\gamma(\cdot)$ is $\widetilde{N}^\gamma(t)=e^{-\gamma t^2}$.\\

\noindent This paper aims at reconstructing the Wigner function $W_\rho$ of a monochromatic light in a cavity prepared in state $\rho$  from $n$ observations. As we cannot measure precisely the quantum state in a single experiment, we perform measurements on $n$ independent identically prepared quantum systems. The measurement carried out on each of the $n$ systems in state $\rho$ is done by QHT as described in Section~\ref{physbackground}. In  practice, the results of such experiments would be $n$  independent identically distributed random variables $(Z_{1}, \Phi_{1}),\dots ,(Z_{n}, \Phi_{n})$ such that
\begin{equation}
         \label{noisy.data}
         Z_\ell:=X_\ell+\sqrt{2\gamma}\,\xi_\ell.
\end{equation}
\noindent with values in $\R\times [0,\pi]$ and distribution $\Pro_{\rho}^\gamma$ having density with respect to the Lebesgue measure on $\R\times [0,\pi]$ equal to $p^\gamma_\rho(\cdot,\cdot)$ defined in (\ref{densitbrui}). For all $\ell=1,\ldots, n$, the $\xi_\ell$'s are independent standard Gaussian random variables, independent of all $(X_{\ell}, \Phi_{\ell})$. \\

\noindent In order to study the theoretical performance of our different procedures, we use the fact that the unknown Wigner function belong to the class of very smooth functions $\Acal(\beta,r,L)$ (similar to those of \cite{Butucea&Guta&Artiles,ABM}) described via its Fourier transform:
\begin{eqnarray}
        \hspace{-1cm}
        \label{class}
         \Acal(\beta,r,L):= \left\{f:\R^2\rightarrow\R,\,\iint |\widetilde{f}(u,v)|^{2}
         e^{2\beta\|(u,v)\|^r} du dv\leqslant(2\pi)^{2}L\right\},
\end{eqnarray}
\noindent where $\widetilde{f}(\cdot,\cdot)$ denotes the Fourier transform with respect to both variables and $\| (u,v)\|=\sqrt{u^2+v^2}$ denote the usual Euclidean scalar norm. Note that this class is reasonable from a physical point of view as the class realistic $\Rcal(C,B,r)$ of density matrix defined in (\ref{eq.classcoeff}) has been translated in terms of Wigner functions by \cite{ABM}. They prove that the fast decay of the elements of the density matrix implies both rapid decay of the Wigner function and of its Fourier transform.

\subsection*{Outline of the results}
\label{Outline.results}

\noindent The problem of reconstructing the quantum state of a light beam has been extensively studied in physical literature 
and in quantum statistics. We mention only papers with theoretical analysis of the performance of their estimation procedure.  Many other physical papers references can be found therein. Methods for reconstructing a quantum state are based on the estimation of either the density matrix $\rho$ or the Wigner function $W_\rho$. In order to compute  the performance of a procedure, a realistic class of quantum states $ \mathcal{R}(C,B,r)$ has defined in many papers as in (\ref{eq.classcoeff}) in which the elements of the density matrix decrease rapidly. From the physical point of view, all the states which have been produced in the laboratory up to date belong to such a class with $r=2$, and a more detailed argument can be found in the paper of \cite{Butucea&Guta&Artiles}.\\

\noindent The estimation of the density matrix from averages of data has been considered in the framework of ideal detection ($\eta=1$ \textit{i.e.} $\gamma=0$) by \cite{Artiles&Gill&Guta} while the noisy setting as investigated by \cite{ABM} for the Frobenius - norm risk.  More recently in the noisy setting, an adaptive estimation  procedure over the classes of quantum states $\mathcal{R}(C,B,r)$, \textit{i.e.} without assuming the knowledge of the regularity parameters, has been proposed by \cite{AMP} and an upper bound for Frobenius - norm risk has been given. The problem of goodness-of-fit testing in quantum statistics has been considered in \cite{MezianiTest}. In this noisy setting, the latter paper derived a testing procedure from a projection-type estimator where the projection is done in $L_2$ distance on some suitably chosen pattern functions.\\

\noindent Note that we may capture some features of the quantum states more easily on the Wigner function $W_\rho$, for instance when this function has significant negative parts, the fact that the quantum state is non classical. \cite{ABM} translate the class $\mathcal{R}(C,B,r)$ in terms of rapid decay of the Fourier transform of its associated Wigner functions as defined in (\ref{class}) by the class $\mathcal A(\beta,r,L)$.  Over this class with $r=1$ and for the problem of pointwise estimation of the Wigner function,  when no noise is present, we mention the work of \cite{Guta&Artiles}. They propose a kernel estimator and derive sharp minimax results over this class.\\

\noindent This paper deals with the problem of reconstruction the Wigner function $W_\rho$ in the context of QHT when taking into account the detection losses occurring in the measurement, leading to an additional Gaussian noise in the measurement data ($\eta\in]1/2, 1]$). The same problem in the noisy setting was treated by \cite{Butucea&Guta&Artiles}, they obtain minimax rates for the pointwise risk over the class $\mathcal A(\beta,r,L)$ for the procedure defined in (\ref{Westimator}).  Moreover, a truncated version of their estimator is proposed by \cite{ABM} where a upper bounds is computed for the $\mathbb{L}_2$ risk over the class $\mathcal A(\beta,r,L)$. The estimation of a quadratic functional of the Wigner function, as an estimator of the purity, was explored in \cite{Meziani}.\\

\noindent The remainder of the article is organized as follows. In Section~\ref{sec.dens.mat}, we establish in Theorem~\ref{bornesupW} the first sup-norm risk upper bound for the estimation procedure (\ref{Westimator}) of the Wigner function while in Theorem~\ref{lower bound} we establish the first minimax lower bounds for the estimation of the Wigner function for the quadratic and the sup-norm risks. These results match our sup-norm upper bounds results up to a logarithmic factor in the sample size $n$.\\

\noindent We propose in Section~\ref{s_adaptation} a Lepski-type procedure that adapts to the unknown smoothness parameters $\beta>0$ and $r\in ]0,2]$ of the Wigner function of interest. The only previous result on adaptation is due to \cite{Butucea&Guta&Artiles} but concerns the simplest case  $r\in ]0,1[$ where the estimation procedure (\ref{Westimator}) with a proper choice of the parameter $h$ independent of $\beta,r$ is naturally minimax adaptive up to a logarithmic factor in the sample size $n$. Theoretical investigations are complemented by numerical experiments reported in Section~\ref{sec.simulations}. The proofs of the main results are defered to the Appendix.

\section{ Wigner function estimation and minimax risk}
\label{sec.dens.mat}

\noindent From now, we work in the practice framework and we assume that $n$ independent identically distributed random pairs $(Z_i,\Phi_i)_{i=1,\ldots,n}$ are observed, where $\Phi_i$ is uniformly distributed in $[0,\pi]$ and the  joint density of $(Z_i,\Phi_i)$ is $p^\gamma_\rho(\cdot,\cdot)$ (see~(\ref{densitbrui})). As \cite{Butucea&Guta&Artiles}, we use the modified the usual tomography kernel in order to take into account the additive noise on the observations and construct a kernel $K_h^\gamma$ which performs both deconvolution and inverse Radon transform on our data, asymptotically such that our estimation procedure is
 \begin{equation}
         \label{Westimator}
          \widehat{W}^\gamma_{h}(q,p) = \frac{1}{2\pi n} \sum_{\ell=1}^n 
          K_h^\gamma \left([ z,\Phi_\ell ] -Z_\ell \right),
\end{equation}
\noindent where $0\leq\gamma<1/4$ is a fixed parameter $h>0$ tends to $0$ when $n \to \infty$ in a proper way to be chosen later. The kernel is defined by
\begin{equation}
           \label{Wnoisy.reg.operator}
            \widetilde K_h^\gamma\left(t \right) =|t|e^{\gamma t^2}\1_{|t|\leq 1/h},
 \end{equation}
\noindent where $z=(q,p)$ and $ [z,\phi]=q \cos \phi + p \sin \phi$. 

\noindent From now,  $\| \cdot \|_{\infty}$ and $\| \cdot \|_2$ and $\| \cdot \|_1$ will denote respectively the sup-norm, the $\mathbb L_2$- norm and the $\mathbb L_1$- norm.  As the sup-norm risk can be trivially bounded as follow
 \begin{eqnarray}
              \hspace{-1cm}
              \label{Bias/Var-decomposition}
                \|\widehat{W}^\gamma_{h}- W_\rho\|_{\infty} \leq \|\widehat{W}^\gamma_{h} -
                \E[\widehat{W}^\gamma_{h}] \|_{\infty}+ \left\|\E\left[\widehat{W}^\gamma_{h}
                \right] -W_\rho\right\|_{\infty},
\end{eqnarray}
\noindent and in order to study the sup-norm risk of our procedure $\widehat{W}^\gamma_{h}$, we study in Proposition~\ref{bias} and \ref{Prop-1-sto-dev}, respectively  the bias term and the stochastic term.

\begin{Proposition}
\label{bias}
\noindent Let $\widehat{W}^\gamma_{h}$ be the estimator of $W_\rho$ defined in (\ref{Westimator}) and $h>0$ tends to $0$ when $n \to \infty$ . Then, 
$$
 \left\|\E\left[\widehat{W}^\gamma_{h}\right] - W_\rho\right\|_{\infty} \leq  \sqrt{\frac{L}{(2\pi)^2\beta r}} h^{(r-2)/2}e^{-\beta h^{-r}} (1+o(1)),
$$
\noindent  where $W_\rho\in\Acal(\beta,r,L)$ defined in (\ref{class}) and $r\in]0,2]$.
\end{Proposition}
The proof is defered to Appendix~\ref{proofbias}.

\begin{Proposition}
\label{Prop-1-sto-dev}
\noindent Let $\widehat{W}^\gamma_{h}$ be the estimator of $W_\rho$ defined in (\ref{Westimator}) and $0<h<1$. Then, there exists a  constant $C_1$, depending only on $\gamma$ such that 
\begin{eqnarray}
         \hspace{-1cm}
         \label{p1}
          \E\left[\|\widehat{W}^\gamma_{h} -\E[\widehat{W}^\gamma_{h}] \|_{\infty}\right]
           \leq C_1 e^{\gamma h^{-2}} \left(\sqrt{\frac{1}{n}} + \frac{1}{n}\right).
\end{eqnarray}
\noindent Moreover, for any $x>0$, we have with probability at least $1-e^{-x}$ that
 \begin{eqnarray}
            \hspace{-1cm}
            \label{p2}
            \|\widehat{W}^\gamma_{h} -\E[\widehat{W}^\gamma_{h}] \|_{\infty} \leq C_2 
            e^{\gamma h^{-2}}\max \left\{ \sqrt{\frac{1+x}{n}}, \frac{1+x}{n}\right\},
\end{eqnarray}
\noindent where $C_2>0$ depends only on $\gamma$.
\end{Proposition}
\noindent The proof is defered to Appendix~\ref{proofProp-1-sto-dev}. The following Theorem  establishes the upper bound of the sup-norm risk. 

\begin{Theorem}
\label{bornesupW}
\noindent Assume that $W_\rho$ belongs to the class $\Acal(\beta,r,L) $ defined in (\ref{class}) for some $r\in ]0,2]$ and $\beta,L >0$. Consider the estimator (\ref{Westimator}) with $h^*=h^{*}(r)$ such that 
 \begin{eqnarray}
          \hspace{-1cm}
          \label{hopti}
          \left\{  \begin{array}{ccc}
           \frac{\gamma}{(h^{*})^2}+ \frac{\beta}{(h^{*})^r}=\frac 12\log(n) & \mbox{ if } 
           &0<r<2,\\
           h^*=\left( \frac{2(\beta+\gamma)}{\log n}\right)^{1/2}  &  \mbox{ if } & r=2.
           \end{array} \right.
 \end{eqnarray} 
\noindent Then we have
$$
 \E\left[ \Vert\widehat{W}^\gamma_{h^*}-W_\rho \Vert_{\infty}\right] \leq Cv_n(r),
 $$
\noindent where $C>0$ can depend only on $\gamma,\beta,r,L$ and the rate of convergence $v_n$ is such that
 \begin{eqnarray}
           \hspace{-1cm}
            \label{psiopti}
            v_n(r)=\left\{ \begin{array}{ccc}
             ( h^*)^{(r-2)/2}e^{-\beta  ( h^*)^{-r} } & \mbox{ if } &0<r<2,\\
             n^{-\frac{\beta}{2(\beta+\gamma)}} &  \mbox{ if } & r=2.
             \end{array} \right.
 \end{eqnarray}  
\end{Theorem}
\noindent Note that for $r\in]0,2)$ the rate of convergence $v_n$ is faster than any logarithmic rate in the sample size but slower than any polynomial rate. For $r=2$, the rate of convergence is polynomial in the sample size.

\noindent \textbf{\textit{Proof of Theorem~\ref{bornesupW}}}: 
\noindent Taking the expectation in (\ref{Bias/Var-decomposition}) and using Propositions~\ref{bias} and \ref{Prop-1-sto-dev}, we get for all $0<h<1$
\begin{eqnarray*}
\hspace{-1cm}
                \E\left[\Vert[\widehat{W}^\gamma_{h} - W_\rho\Vert_{\infty}\right] 
                \hspace{-0.15cm}&\leq&\hspace{-0.15cm}
                 \E\left[\Vert\widehat{W}^ \gamma_{h} -\E[\widehat{W}^\gamma_{h}] 
                 \Vert_{\infty}\right] +\Vert\E\left[\widehat{W}^\gamma_{h}\right] - W_\rho
                 \Vert_{\infty}\\
                 \hspace{-0.15cm}&\leq&\hspace{-0.15cm}
                 C e^{\gamma h^{-2}} \sqrt{\frac{1}{n}} (1+o(1))+ C_B h^{(r-2)/2}
                 e^{-\beta h^{-r}} (1+o(1))
\end{eqnarray*}
\noindent where $C_B=\sqrt{\frac{L}{(2\pi)^2\beta r}}$, $h\rightarrow 0$ as $n\rightarrow \infty$ and $W_\rho\in\Acal(\beta,r,L) $. The optimal bandwidth parameter $ h^*(r):= h^*$ is such that
 \begin{eqnarray}
         \hspace{-1cm}
        \label{tradeoff}
         h^*=\arg \inf_{h>0}\left\{C_B h^{(r-2)/2}e^{-\beta h^{-r}}+ 
          C e^{\gamma h^{-2}} \sqrt{\frac{1}{n}}  \right \}.
 \end{eqnarray}
 \noindent Therefore, by taking derivative, we get
 $$
  \frac{\gamma}{ ( h^*)^2}+ \frac{\beta}{ ( h^*)^r}=\frac 12\log(n)+ C(1+o(1)).
  $$
 \noindent  By plugging  the result in (\ref{tradeoff}) for $0<r<2$ we have
 $$
  ( h^*)^{(r-2)/2}e^{-\beta  ( h^*)^{-r}}=( h^*)^{(r-2)/2}\frac{1}{\sqrt{n}}e^{\gamma ( h^*)^{-2}}.
  $$
 \noindent It comes that the bias term is much larger than the stochastic term for $0<r<2$. It is easy to see that for $r=2$, we have  $h^*=\left( \frac{2(\beta+\gamma)}{\log n}\right)^{1/2} $ and that the the bias term and the stochastic term are of the same order. $\square$\\
 
 \noindent  We derive now a minimax lower bound. We consider specifically the case $r=2$ since it is relevant with quantum physic applications, but our results can easily be generalized to the case $r\in]0,2]$. However, similar arguments can be applied to the case $0<r<2$. The only known lower bound result for the estimation of a Wigner function is due to \cite{Butucea&Guta&Artiles} and concerns the pointwise risk. In Theorem~\ref{lower bound} below, we obtain the first minimax lower bounds for the estimation of a Wigner function $W_{\rho} \in \Acal(\beta,2,L)$ with the quadratic and sup-norm risks.

\begin{Theorem}\label{lower bound}
Assume that $(Z_1,\Phi_1),\cdots,(Z_n,\Phi_n)$ coming from the model (\ref{Z}) with $\gamma \in [0,1/4[$. Then, for any $\beta,L>0$ and $r=2$ there exists a constant $c:=c(\beta,L,\gamma)>0$ such that for $n$ large enough
\begin{eqnarray*}\hspace{-1cm}
                    \inf_{\widehat W_n}\sup_{W_{\rho} \in \Acal(\beta,2,L)} 
                    \E \| \widehat W_n - W_{\rho} \|_p^2
                    \hspace{-0.15cm}&\geq&\hspace{-0.15cm}
                    \left\{  \begin{array}{ccc}
                    c n^{-\frac{\beta}{2(\beta+\gamma)}}\log^{-3/2}(n) & \mbox{ if } &p=\infty,\\
                     c n^{-\frac{\beta}{\beta+\gamma}}&  \mbox{ if } &p=2.
                      \end{array} \right.\\
\end{eqnarray*}
\noindent where the infimum is taken over all possible estimators $\widehat W_n$ based on the i.i.d. sample $\left\{(Z_i,\Phi_i)\right\}_{i=1}^n$.
\end{Theorem} 
\noindent The proof is defered to Appendix~\ref{proofLB}. This theorem guarantees that the sup-norm upper bound derived in Theorem \ref{bornesupW}  and the quadratic risk upper bound in the paper of \cite{ABM} are minimax optimal up to a logarithmic factor in the sample size. We believe that the logarithmic factors for both cases are artefact of the proofs.

\section{Adaptation to the smoothness}
\label{s_adaptation}

\noindent As we see in (\ref{tradeoff}), the optimal choice of the bandwidth $h^*$ depends on the unknown smoothness $\beta$.  For any $0<r\leq 2$, we propose here to implement a Lepski type procedure to select an adaptive bandwidth $h$. We will show that the estimator obtained with this bandwidth achieves the optimal minimax rate up to a logarithmic factor. Our adaptive procedure is implemented in Section~\ref{sec.simulations}.\\

\noindent Let $M\geq 2$, and $0< h_M < \cdots < h_1< 1$ a grid of $]0,1[$, we build  estimators $\widehat{W}^\gamma_{h_m}$ associated to bandwidth $h_m$ for any $1\leq m \leq M$. For any fixed $x>0$, let us define  
$r_n(x) = \max\left( \sqrt{\frac{1+x}{n}} , \frac{1+x}{n}\right)$. 
We denote by $\Lcal_\kappa(\cdot)$, the Lepski functional such that
\begin{eqnarray}
             \hspace{-1cm}
             \label{Lep-proc1}
             \mathcal L_{\kappa}(m)\hspace{-0.15cm}&=&\hspace{-0.15cm}max_{j>m}
             \left\{  \|\widehat{W}^\gamma_{h_m} - \widehat{W}^\gamma_{h_j}\|_{\infty} - 
             2 \kappa e^{\gamma h_j^{-2}} r_n(x+\log M) \right\}\nonumber\\
             && + 2 \kappa e^{\gamma h_m^{-2}}r_n(x+\log M),
\end{eqnarray}
\noindent where $\kappa>0$ is a fixed constant. Therefore, our final adaptive estimator denoted by $\widehat{W}^\gamma_{h_{\mathbf{\widehat m }}}$ will be the estimator defined in (\ref{Westimator}) for the bandwidth $h_{\mathbf{\widehat m }}$. The bandwidth $h_{\mathbf{\widehat m }}$ is such that 
\begin{eqnarray}
     \label{Lep-proc2}
     \mathbf{\widehat m }= \mathrm{argmin}_{1\leq m \leq M}\mathcal L_{\kappa}(m).
\end{eqnarray}

\begin{Theorem}\label{thm-adaptation}
Assume that $W_\rho \in \Acal(\beta,r,L)$. Take $\kappa >0$ sufficiently large and $M\geq 2$. Choose $0<h_M < \cdots < h_1< 1$. Then, for the bandwidth $h_{\mathbf{\widehat m }}$ with $\mathbf{\widehat m }$ defined in (\ref{Lep-proc2}) and  for any $ x >0$, we have with probability at least $1-e^{-x}$
\begin{eqnarray}
         \hspace{-1cm}
         \label{p3}
         \|\widehat{W}^\gamma_{h_{\mathbf{\widehat m }}}  - W_{\rho}\|_{\infty}\leq 
         C\min_{1\leq m \leq M}\left\{  h_m^{r/2-1} e^{-\frac{\beta}{h_m^{r}}} +
           e^{\gamma h_{m}^{-2}}r_n(x+\log M) \right\},
\end{eqnarray}
where $C>0$ is a constant depending only on $\gamma,\beta,r,L$.\\
In addition, we have in expectation
\begin{eqnarray}
         \hspace{-1cm}
         \label{p4}
         \E\left[\|\widehat{W}^\gamma_{h_{\mathbf{\widehat m }}}  - W_{\rho}\|_{\infty}
         \right] \leq C'\min_{1\leq m \leq M}\left\{  h_m^{r/2-1} e^{-\frac{\beta}{h_m^{r}}} + 
          e^{\gamma h_{m}^{-2}}r_n(\log M) \right\},
\end{eqnarray}
where $C'>0$ is a constant depending only on $\gamma,r,\beta,L$.
\end{Theorem}
\noindent The proof is defered to the Appendix~\ref{proofthm-adaptation}. \\

\noindent  The idea is now to build a sufficiently fine grid $0<h_M< \cdots<h_1<1$ to achieve the optimal rate of convergence simultaneously over $r\in (0,2]$ and $\beta >0$. Take $M=\lfloor  \sqrt{\log n/(2\gamma)}  \rfloor$. We consider the following grid for the bandwitdh parameter $h$:
\begin{eqnarray}
           \hspace{-1cm}\label{h-grid}
           h_1 =  1/ 2  ,\quad h_m =\frac {1}{ 2}\left( 1 - (m-1)\sqrt{\frac{2\gamma}{\log n}}
           \right),\quad 1\leq m \leq M.
\end{eqnarray}
We build the corresponding estimators $\widehat W^{\gamma}_{h_m}$ and we apply the Lepski procedure (\ref{Lep-proc1})-(\ref{Lep-proc2}) to obtain the estimator $ \widehat W^{\gamma}_{h_{\mathbf{\widehat m }}}$. The next result guarantees that this estimator is minimax adaptive over the class 
$$
\Omega:=\left\{  (\beta,r,L),\; \beta >0, \,  0<r \leq 2,\; L>0  \right\}.
$$

\begin{Corollary}
\label{corr1}
Let the conditions Theorem \ref{thm-adaptation} be satisfied. Then the estimator $\widehat W^{\gamma}_{h_{\mathbf{\widehat m }}}$ for the bandwidth $h_{\mathbf{\widehat m }}$ with $\mathbf{\widehat m }$ defined in (\ref{Lep-proc2})
 and for any $(\beta,r,L)\in\Omega$ satisfies 
$$
\mathrm{limsup}_{n\rightarrow\infty} \sup_{W_\rho\in\Acal(\beta,r,L)}\E\left[ \Vert \widehat W^{\gamma}_{h_{\mathbf{\widehat m }}}-W_\rho\Vert_\infty \right]\leq Cv_n(r),
$$
\noindent where  $v_n(r)$  is the rate defined in (\ref{psiopti}) and $C$ is a positive constant depending only on $r$, $L$, $\beta$ and $\gamma$.
\end{Corollary}

\noindent\textbf{ \textit{Proof of Corollary~\ref{corr1}} }: 
\noindent First note that for all $m=1,\cdots,M$ and as 
$$
h_m \in](\gamma/(2\log n))^{1/2},1/ 2],
$$
 the bias term $ h_m^{r/2-1} e^{-\frac{\beta}{h_m^{r}}} $ is larger than the stochastic term $e^{\gamma h_{m}^{-2}}r_n(\log M)$ up to a numerical constant. Let define  
$$
\widetilde{m}:=\arg\max_{1\leq m\leq M}\{|h_m-h^*|: h_m\leq h^*\},
$$
\noindent where $\widetilde m$ is well defined as 
\begin{eqnarray*}
\hspace{-1cm}
                 \frac{h_M}{h^*}\hspace{-0.15cm}&=&\hspace{-0.15cm}\frac{(1/ 2)
                 \left(1-M(2\gamma/\log n)^{1/2} +(2\gamma/\log n)^{1/2}\right)} 
                 {\left(  \log n/(2\gamma) -  (\beta/\gamma) (h^*)^{-r}\right )^{-1/2}  }\\
                 \hspace{-0.15cm}&=&\hspace{-0.15cm}\frac{1  }{ 2  }\left(1-M 
                 +\left((\log n)/(2\gamma)\right) ^{1/2}\right)\left(  1 -  (2\beta/(\log(n))
                  (h^*)^{-r}\right )^{1/2}. 
\end{eqnarray*}
\noindent Moreover, as $0\leq\left((\log n)/(2\gamma)\right) ^{1/2} -M\leq 1$ we get 
$$
\frac{h_M}{h^*}\leq\left(  1 -  (2\beta/(\log(n)) (h^*)^{-r}\right )^{1/2} \leq 1.
$$
\noindent Therefore, from (\ref{p4}), 
\begin{eqnarray*}
\hspace{-1cm}
                    \E\left[\|\widehat{W}^\gamma_{h_{\mathbf{\widehat m }}}  - 
                    W_{\rho}\|_{\infty}\right] \hspace{-0.15cm}&\leq&\hspace{-0.15cm} 
                    C  h_{\widetilde{m}}^{r/2-1} e^{-\frac{\beta}{h_{\widetilde{m}}^{r}}}
                    \leq C  h_{\widetilde{m}}^{r/2-1} e^{-\frac{\beta}{h_{\widetilde{m}}^{r}}}  
                    v_n(r)v_n(r)^{-1}\\
                    \hspace{-0.15cm}&=&\hspace{-0.15cm} C \left(\frac{  h_{\widetilde{m}}}
                    {h^*}\right)^{r/2-1}e^{-\beta (h_{\widetilde{m}}^{-r}  - (h^*)^{-r}) }  v_n(r).
\end{eqnarray*}
\noindent By the definition of $\widetilde m$, it comes that 
$h_{\widetilde{m}}^{-r}   >(h^*)^{-r}$, then 
 \begin{eqnarray*}
 \hspace{-1cm}
                    \E\left[\|\widehat{W}^\gamma_{h_{\mathbf{\widehat m }}}  -
                     W_{\rho}\|_{\infty}\right] \hspace{-0.15cm}&\leq&\hspace{-0.15cm} 
                      C \left(\frac{  h_{\widetilde{m}}}{h^*}\right)^{r/2-1} v_n(r)= 
                      C \left(\frac{  h_{\widetilde{m}} -h^* }{h^*} +1\right)^{r/2-1} v_n(r).
\end{eqnarray*}
\noindent By construction $|h_{\widetilde{m}} -h^*|\leq (\gamma/(2\log n))^{1/2}$, then we have
\begin{eqnarray*}
\hspace{-1cm}
                      \E\left[\|\widehat{W}^\gamma_{h_{\mathbf{\widehat m }}}  - 
                      W_{\rho}\|_{\infty}\right] \hspace{-0.15cm}&\leq&\hspace{-0.15cm} 
                      C \left(1-\frac{ (\gamma/(2\log n))^{1/2}   }{h^*} \right)^{r/2-1} v_n(r).
\end{eqnarray*}
\noindent  As $(h^*)^{-1}\leq \left(  \log n/(2\gamma) \right)^{1/2}$,  it holds
$1-\frac{ (\gamma/(2\log n))^{1/2}   }{h^*} \geq 1/2.$ Therefore as  $r/2-1<0$, the result follow
\begin{center}
$\quad\quad \quad \quad \quad \quad \quad\quad \quad \quad\E\left[\|\widehat{W}^\gamma_{h_{\mathbf{\widehat m }}}  - W_{\rho}\|_{\infty}\right] \leq C v_n(r).\quad \quad \quad \quad \quad \quad \quad \quad \quad \quad \quad \quad \quad \quad \quad \quad \   \square$
\end{center}

\section{Experimental evaluation}
\label{sec.simulations}

\newcommand{\EstimLepski}{\widehat{W}^\gamma_{h_{\mathbf{\widehat m}}}}
\newcommand{\Estim}{\widehat{W}^\gamma_{h}}

We test our method on two examples of Wigner functions, corresponding to the single-photon and the Schrödinger's cat states, and that are respectively  defined as
\begin{align*}
	W_\rho(q,p) &= - (1-2 (q^2+p^2) ) e^{-q^2-p^2},  \\ 
	W_\rho(q,p) &=  \frac{1}{2} e^{-(q-q_0)^2-p^2}   + 
				\frac{1}{2} e^{-(q+q_0)^2-p^2} + 
              \cos(2 q_0 p) e^{-q^2-p^2}.
\end{align*}
\noindent We used $q_0=3$ in our numerical tests. The toolbox to reproduce the numerical results of this article is available online\footnote{\url{https://github.com/gpeyre/2015-AOS-AdaptiveWigner}}. Following the paper of~\cite{Butucea&Guta&Artiles} and in order to obtain a fast numerical procedure, we implemented the estimator $\Estim$ defined in~\eqref{Westimator} on a regular grid. More precisely, 2-D functions such as $W_\rho$ are discretized on a fine 2-D grid of $256 \times 256$ points. We use the Fast Slant Stack Radon transform of~\cite{SlantStack}, which is both fast and faithful to the continuous Radon transform $\mathcal{R}$. It also implements a fast pseudo-inverse which accounts for the filtered back projection formula~\eqref{Westimator}. The filtering against the 1-D kernel~\eqref{Wnoisy.reg.operator} is computed along the radial rays in the Radon domain using Fast Fourier transforms. We computed the Lepski functional~\eqref{Lep-proc1} using the values $x=\log(M)$ and $\kappa = 1$.

\newpage

\newcommand{\FigureTable}[2]{ 
\begin{tabular}{@{}c@{\hspace{2mm}}c@{}c@{}}
	\FigErrPlots{#1}{#2}{error} & 
	\FigPlotsLS{#1}{}{}	&
	\FigPlotsSURF{#1}{}{} \\
	$\|\Estim-W_\rho\|_{\infty}/\|W_\rho\|_{\infty}$ as a function of $h$ & 
	$W_\rho$ (2-D display) &
	$W_\rho$ (3-D display) \\
	\FigErrPlots{#1}{#2}{lepski-histo} & 
	\FigPlotsLS{#1}{#2}{-lepski-} &
	\FigPlotsSURF{#1}{#2}{-lepski-} \\
	Histogram of the repartition of $h_{\mathbf{\widehat m}}$ & 
	$\EstimLepski$ (2-D display) &
	$\EstimLepski$ (3-D display)  \\
	\end{tabular}
}

\begin{figure}[h!]
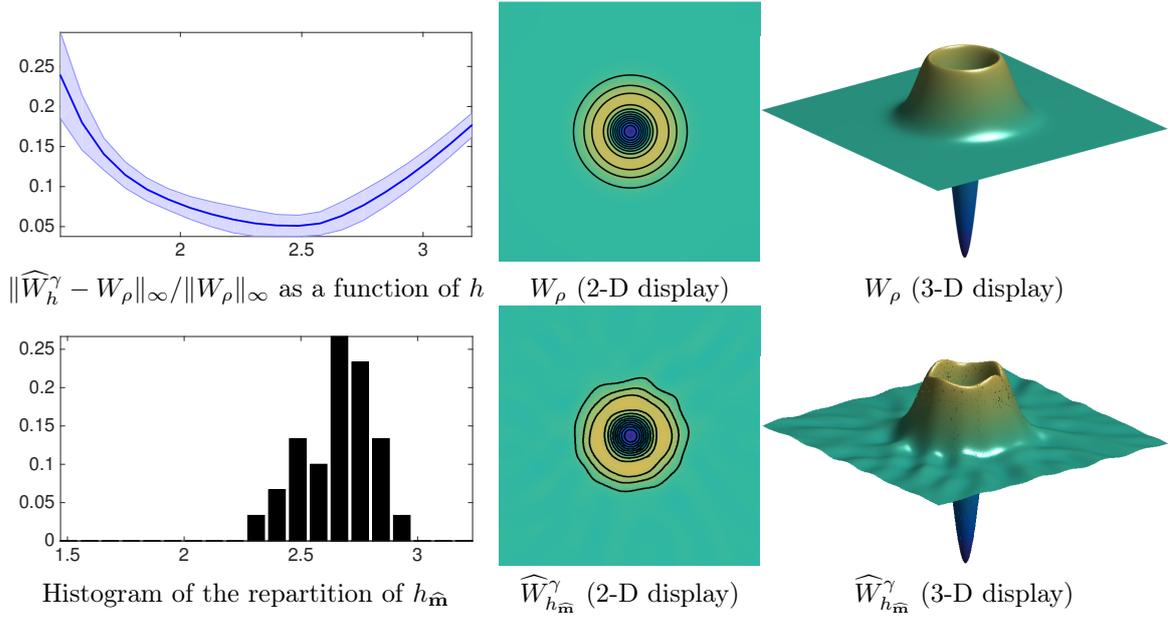

	\centering
	\FigureTable{single-photon}{eta90-n100e3}
	\caption{%
		\small{Single photon cat state estimation, with $\eta=0.9$, $n=100 \times 10^3$.  
		\textit{Left, top:} display of $\|\Estim-W_\rho\|_{\infty}/\|W_\rho\|_{\infty}$ as a function of $1/h$.
			The central curve is the mean of this quantity, while
			the shaded area displays the $\pm 2 \times $ standard deviation of this quantity. 
		\textit{Left, right:} histogram of the empirical repartition of $\mathbf{\widehat m}$
		computed by the Lepski procedure~\eqref{Lep-proc2}. 
		\textit{Center:} display as a 2-D image using level sets of $W_\rho$ (top) 
		and $\EstimLepski$ (bottom).
		\textit{Right:} same, but displayed as an elevation surface.}
	}
   \label{fig-single-photon}
\end{figure}
\begin{figure}[h!]
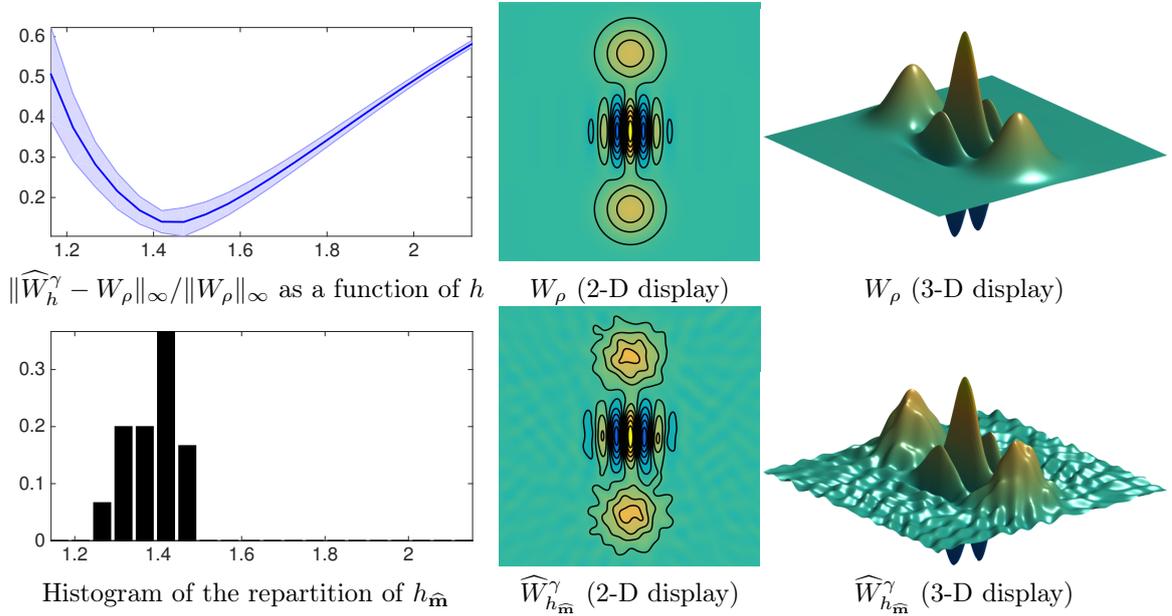

	\centering
	\FigureTable{schrodinger-cat}{eta90-n500e3}
	\caption{%
		\small{Schrödinger's cat state estimation, with $\eta=0.9$, $n=500 \times 10^3$.  
		We refer to Figure~\ref{fig-single-photon} for the description of the plots.} 
	}
   \label{fig-schrodinger-cat}
\end{figure}

\newpage

\noindent Figures~\ref{fig-single-photon} and~\ref{fig-schrodinger-cat} reports the numerical results of our method on both test cases.  The left part compares the error $\|\Estim-W_\rho\|_{\infty}$ (displayed as a function of $h$) to the parameters $h_{\mathbf{\widehat m}}$ selected by the Lepski procedure~\eqref{Lep-proc2} . The error $\|\Estim-W_\rho\|_{\infty}$ (its empirical mean and its standard deviation) is computed in an ``oracle'' manner (since for these examples, the Wigner function to estimate $W_\rho$ is known) using 20 realizations of the sampling for each tested value $(h_i)_{i=1}^M$. The histogram of values $h_{\mathbf{\widehat m}}$ is computed by solving~\eqref{Lep-proc1} for 20 realizations of the sampling.  This comparison shows, on both test cases, that the method is able to select a parameter value $h_{\mathbf{\widehat m}}$ which lies around the optimal parameter value (as indicated by the minimum of the $L^\infty$ error). The central and right parts show graphical displays of $\EstimLepski$, where $\mathbf{\widehat m }$ is selected using the Lepski procedure~\eqref{Lep-proc2}, for a given sampling realization.


\appendix
\section{Proof of Propositions}
\label{proofbias}

\subsection{Proof of Proposition~1}
\label{proofbias}

\noindent First remark that by the Fourier transform formula for $w=(q,p)\in\R^2$ and $x=(x_1,x_2)$
 \begin{eqnarray}
\label{bias1}
               W_\rho(w)\hspace{-0.15cm}&=&\hspace{-0.15cm}
               \frac{1}{(2\pi)^2}\iint \Fourier{W}_\rho(x)e^{-i(qx_1+px_2)}dx.
 \end{eqnarray}
\noindent Let $\widehat{W}^\gamma_{h}$ be the estimator of $W_\rho$ defined in (\ref{Westimator}), then    
\begin{eqnarray*}
\hspace{-1cm}
               \E\left[\widehat{W}^\gamma_{h}(w)\right]\hspace{-0.15cm}&=&\hspace{-0.15cm}
               \frac{1}{2\pi} \E\left[K_{h}^{\gamma}([w, \Phi_1] - Z_1)\right]
             =\frac{1}{2\pi}\int_{0}^{\pi}\hspace{-0.25cm}\int K_{h}^{\gamma}([w,\phi] - z)
             p_\rho^\gamma(z,\phi)dzd\phi\\
             \hspace{-0.15cm}&=&\hspace{-0.15cm}\frac{1}{2\pi} \int_0^\pi  K_{h}^{\gamma}
             *p_\rho^\gamma(\cdot,\phi)([w,\phi]) d\phi.
\end{eqnarray*}
\noindent In the fourier domain, the convolution becomes a product, combining with (\ref{fourierproun}), we obtain
\begin{eqnarray*}
\hspace{-1cm}
                  \E\left[\widehat{W}^\gamma_{h}(w)\right]\hspace{-0.15cm}&=&
                  \hspace{-0.15cm}\int_0^\pi \frac{1}{(2\pi)^2}\int
                  \Fourier{K}_{h}^{\gamma}(t) \Fcal_1[p_\rho^\gamma(\cdot,\phi)](t)e^{-it[w, \phi]}
                  dt d\phi.
\end{eqnarray*}
 \noindent As $\widetilde{N}^\gamma(t)=e^{-\gamma t^2}$, the definition (\ref{Wnoisy.reg.operator}) of the kernel combining with (\ref{fourierproun}) gives 
\begin{eqnarray*}\hspace{-1cm}
             \E\left[\widehat{W}^\gamma_{h}(w)\right]  \hspace{-0.15cm}&=&\hspace{-0.15cm}
             \int_0^\pi \frac{1}{(2\pi)^2}\int
          \Fourier{K}_{h}^{\gamma}(t)\Fourier{W}_\rho(t\cos(\phi),t\sin(\phi))
          \widetilde{N}^\gamma(t)e^{-it[w,\phi]}dtd\phi\\
         \hspace{-0.15cm}&=&\hspace{-0.15cm}\int_0^\pi \frac{1}{(2\pi)^2}\int_{ |t|\leq 1/h} |t|
         \Fourier{W}_\rho(t\cos(\phi),t\sin(\phi))e^{-it[w,\phi]}dtd\phi.
 \end{eqnarray*}
\noindent  Therefore, by the change of variable $x=(t\cos(\phi),t\sin(\phi))$, it comes
\begin{eqnarray}
       \label{bias2}
         \E\left[\widehat{W}^\gamma_{h}(w)\right]  \hspace{-0.15cm}&=&\hspace{-0.15cm}
         \frac{1}{(2\pi)^2}\int_{ ||x||_2\leq 1/h}
        \Fourier{W}_\rho(x)e^{-i(qx_1+px_2)}dx.
\end{eqnarray}
\noindent  From equations (\ref{bias1}) and (\ref{bias2}), we have
\begin{eqnarray*}
 \hspace{-0.15cm}
             \left|\E\left[\widehat{W}^\gamma_{h}(w)\right]-W_\rho(w) \right|\hspace{-0.15cm}
             &\leq&\hspace{-0.15cm}\frac{1} {(2\pi)^2} \int_{ ||x||_2>1/h}
             \left|  \Fourier{W}_\rho(x)\right |dx \\
            \hspace{-0.15cm}&\leq&\hspace{-0.15cm} \frac{1}{(2\pi)^2}
            \left[\iint \left|  \Fourier{W}_\rho(x)\right |^2 e^{2\beta  ||x||_2^{r/2}}dx\right]^{1/2}
             \left[\int_{ ||x||_2>1/h}e^{-2\beta ||x||_2^{r/2}}dx \right]^{1/2}\\
              \hspace{-0.15cm}&\leq&\hspace{-0.15cm} \sqrt{\frac{L}{(2\pi)^2\beta r}}
               h^{(r-2)/2}e^{-\beta h^{-r}} (1+o(1)),\quad h\rightarrow 0
\end{eqnarray*}
\noindent as $W_\rho\in\Acal(\beta,r,L)$ the class defined in (\ref{class}).

\subsection{{Proof of Proposition~2}}
\label{proofProp-1-sto-dev}

\noindent The following Lemma is needed to prove the Proposition~\ref{Prop-1-sto-dev}.
%
\begin{Lemma}
\label{VC-class}
        Let $\delta_h:= h^{-1}e^{\frac{\gamma}{h^2}}>0$ for any $0<h\leq 1$, then the class 
\begin{eqnarray}
              \label{un}
              \mathcal H_h\hspace{-0.15cm}&=&\hspace{-0.15cm} \{ \delta_h^{-1} K_{h}^{\eta}(\cdot - t), t\in \R\},\quad h>0
\end{eqnarray}
\noindent is uniformly bounded by $U:=\frac{h}{2 \gamma \pi}$.  Moreover, for every $0 < \epsilon<A$ and for finite positive constants $A,v$ depending only on $\gamma$, 
\begin{eqnarray}
       \label{deux}
       \sup_{Q} N(\epsilon,\mathcal H_h, L^2(Q)) \hspace{-0.15cm}&\leq&\hspace{-0.15cm} (A/\epsilon)^v,
\end{eqnarray}
\noindent  where the supremum extends over all probability measures $Q$ on $\R$.
\end{Lemma}
\noindent The proof of this Lemma can be found in~\ref{proofVC-class}. To prove  (\ref{p1}), we have to bound the following quantity :
\begin{eqnarray}
       \hspace{-1cm}\E[|K_{h}^{\gamma}([z,\phi] - Y/\sqrt{\eta})|^2] \hspace{-0.15cm}&\leq&\hspace{-0.15cm}
       \|K^{\gamma}_{h}\|_{\infty}^2 \leq \|\widetilde K_{h}^{\gamma}\|_1^2
         = \left[  \int_{|t|\leq h^{-1}} |t|e^{\gamma t^2} dt  \right]^2 \nonumber\\
          \label{a2}
        \hspace{-0.15cm}&=&\hspace{-0.15cm} \left[ 2 \int_{0}^{h^{-1}} te^ {\gamma t^2} dt  \right]^2
       = \left( \gamma^{-1} e^{\gamma h^{-2}} -\frac{1}{\gamma}\right)^2 \leq \frac{1}
        {\gamma^2} e^{2\gamma h^{-2}}.
\end{eqnarray}

\noindent Moreover for $\delta_h=h^{-1}e^{\gamma h^{-2}}$, we have
\begin{eqnarray}
\label{sigma}
\delta^{-2}_h\E[|K_{h}^{\eta}([z,\phi] - Y/\sqrt{\eta})|^2] \hspace{-0.15cm}&\leq&\hspace{-0.15cm}\frac{h^{2}}{\gamma^2}.
\end{eqnarray}

\noindent By Lemma~\ref{VC-class}, it comes that the class $\mathcal H_h$ is VC. Hence, we can apply (57) in the paper of \cite{GN2010-aop} to get
\begin{eqnarray}
         \hspace{-0.5cm}
           \E\left[\|\widehat{W}^\gamma_{h} -\E[\widehat{W}^\gamma_{h}] \|_{\infty}\right] \hspace{-0.15cm}&=&\hspace{-0.15cm}
            \E\sup_{z\in \R^2} \left|\frac{1}{2\pi n}\sum_{l=1}^n K_{h}^{^\gamma}\left([z,\phi_l]-
            Z_\ell\right) - \E\left[K_{h}^{\gamma} ([z,\phi_\ell]- Z_\ell)\right] \right|\nonumber \\
            \hspace{-0.15cm}&=&\hspace{-0.15cm}\frac{\delta_h}{2\pi n} \E\sup_{z\in \R^2} \left|\sum_{l=1}^n \delta_h^{-1} 
            K_{h}^{^\gamma}([z,\phi_\ell]- Z_\ell)- \E[\delta_h^{-1} K_{h}^{^\gamma}([z,\phi_\ell]-
             Z_\ell)] \right|\nonumber \\
             \label{a1}
             \hspace{-0.15cm}&\leq&\hspace{-0.15cm}\frac{C(\gamma)\delta_h}{2\pi n} \left(  \sigma \sqrt{n \log \frac{AU}{\sigma}} 
             + U\log \frac{AU}{\sigma}\right),
\end{eqnarray}
where $U =\frac{ h}{2 \gamma\pi} $  is the envelop of the class $\mathcal H_h$ defined in Lemma~\ref{VC-class}. By choosing 
      $$
      \sigma^2 : = \frac{h}{\gamma}\geq \sup_{z\in \R^2} \E\left[\left(\delta_h^{-1}K_h^{\eta}
      ([z,\phi_l]- \frac{Y_l}{\sqrt{\eta}})\right)^2\right]
      $$
in (\ref{a1}) we get the result in expectation (\ref{p1}). \\

\noindent To prove the result in probability (\ref{p2}), we use Talagrand's inequality as in Theorem 2.3 of \cite{Bousquet2002}. Let us define
$$Z := \frac{n \gamma}{h \delta_h} \|\widehat{W}^\gamma_n -\E[\widehat{W}^\gamma_n] \|_{\infty}.$$

\noindent In view of the previous display (\ref{a2}), we have
\begin{eqnarray*}
       \hspace{-1.3cm}\mathrm{Var}\left(\gamma (h\delta_h)^{-1} \left\vert K_h^{\gamma}(\cdot - 
       \frac{Y_1}{\sqrt{\gamma}} ) - \mathbb E\left [ K_h^{\gamma}(\cdot -\frac{Y_1}{\sqrt
       {\gamma}})\right]\right\vert\right)
       \hspace{-0.15cm}&\leq&\hspace{-0.15cm} \gamma^2 (h\delta_h)^{-2}\E\left[\left\vert   K_h^{\gamma}(\cdot - \frac{Y_1}{\sqrt
       {\gamma}} )    \right\vert^2\right]\nonumber\\
       \label{a3}
       \hspace{-0.15cm}&\leq&\hspace{-0.15cm} \gamma^2 (h\delta_h)^{-2} \frac{1}{\gamma^2} e^{2\gamma h^{-2}} =1.
\end{eqnarray*}
\noindent As $U=\frac{h}{2 \gamma \pi}$ and by (\ref{U}), it comes
\begin{eqnarray*}
       \hspace{-1.3cm}
\gamma (h\delta_h)^{-1}\left \Vert K_h^{\gamma} (\cdot - \frac{Y_1}{\sqrt{\gamma}} )- \mathbb E [ K_h^{\gamma}(\cdot - \frac{Y_1}{\sqrt{\gamma}} )]\right \Vert_{\infty}\leq  \gamma (h\delta_h)^{-1} \Vert K_h^{\gamma} \Vert_\infty  \leq   \gamma h^{-1} U\leq    1.
\end{eqnarray*}
\noindent Then, for  any $x>0$ and with probability at least $1-e^{-x}$, we obtain
\begin{eqnarray*}\hspace{-1.5cm}
        Z \hspace{-0.15cm}&\leq&\hspace{-0.15cm}\E\left[ Z\right] + \sqrt{2xn + 4 x\E[Z]} + \frac{x}{3}
       \leq \E\left[ Z\right] + \sqrt{2xn} + 2\sqrt{x\E[Z]} + \frac{x}{3}\\
       \hspace{-0.15cm}&\leq&\hspace{-0.15cm} 2\E\left[ Z\right]  + \sqrt{2xn}  + \frac{4x}{3},
\end{eqnarray*}
\noindent where we have used the decoupling inequality $2ab \leq a^2 + b^2$ with $a=\sqrt{x}$ and $b=\sqrt{\E[Z]}$. Thus, with probability at least $1-e^{-x}$, we get
\begin{eqnarray*}\hspace{-1cm}
             \|\widehat{W}^\gamma_{h} -\E[\widehat{W}^\gamma_{h}] \|_{\infty}=  \frac{h \delta_h}{n 
             \gamma} Z&\leq 2\E\left[ \|\widehat{W}^\gamma_{h} -\E[\widehat{W}^
             \gamma_{h}] \|_{\infty}\right] + \frac{e^{\gamma h^{-2}}}{\gamma}\left(\sqrt{2\frac{x}{n}}
              + \frac{4x }{3 n}\right). 
\end{eqnarray*}
Plugging our control (\ref{p1}) on $\E[\|\widehat{W}^\gamma_{h} -\E[\widehat{W}^\gamma_{h}] \|_{\infty}]$, the result in probability follows.

\section{Proof of Theorem~2 - Lower bounds}
\label{proofLB}

\subsection{Proof of Theorem~2 - Lower bounds for the $\mathbb L_2$-norm}
\label{proofLB2}

The proof for the minimax lower bounds follows a standard scheme for deconvolution problem as in the paper of \cite{Butucea&Guta&Artiles,LN2011}. However, additional technicalities arise to build a proper set of Wigner functions and then to derive a lower bound. From now on, for the sake of brevity, we will denote $\Acal(\beta,2,L)$ by $\Acal(\beta,L)$ as we consider the practice case $r=2$. Let $W_0\in\Acal(\beta,L)$ be a Wigner function. Its associated density function will be denoted by $p_0(x,\phi)=\frac 1\pi \Rcal[W_0](x,\phi)\1_{[0,\pi]}(\phi)$.\\

\noindent Let $M=\lfloor\sqrt{\log n}\rfloor$ be the integer part of $\log n$, and 
\begin{eqnarray}
\label{delta}
\delta:=\log^{-1}( n).
\end{eqnarray}
\noindent We suggest the construction of a family of $M$ Wigner functions such that for all $m=1,\cdots,M$ and $w\in\R^2$:
$$
W_{m,h}(w) = W_{0}(w) + V_{m,h}(w),\quad 1\leq m \leq M,
$$
depending on a parameter $h = h(n)\rightarrow 0$ as $n\rightarrow \infty$. The construction of $W_0$ and $V_{m,h}$ are discussed in Appendix~\ref{W0} and \ref{blabla}.  We denote by 
$$p_{m,h}(x,\phi)=\frac 1\pi \Rcal[W_{m,h}](x,\phi)\1_{[0,\pi]}(\phi)$$
the associated density function of the Wigner function $W_{m,h}$. As we consider the noisy framework (\ref{Z}) and in view of (\ref{densitbrui}), we set for all $1\leq m \leq M$
\begin{center}
$p_{m,h}^{\gamma}(z,\phi) = \left[ p_{m,h}(\cdot,\phi) \ast N^{\gamma} \right](z)\quad$ and $\quad p_{0}^{\gamma}(z,\phi) = \left[ p_{0}(\cdot,\phi) \ast N^{\gamma} \right](z)$.
\end{center}
  If the following conditions $\mathbf{(C1)}$ to $\mathbf{(C3)}$ are satisfied, then Theorem 2.6 in the book of \cite{Tsybakovlivre} gives the lower bound.

\begin{itemize}
\item[$\mathbf{(C1)}$]  For all $m=1\cdots M $, $\quad W_{m,h}\in\Acal(\beta,L).$
\item[$\mathbf{(C2)}$]   For any $1\leq k\neq m\leq M,$ we have for  $ ||W_{k,h}-W_{m,h}||^2_2\geq 4\varphi^2_{n},$ with  $\varphi^2_{n}=\mathcal O\left(n^{-\frac{\beta}{\beta+\gamma}}\right)$.
\item[$\mathbf{(C3)}$] For all $ 1\leq m \leq M$,
\begin{eqnarray*}
       \hspace{-1cm}
n\mathcal X^2(p_{m,h}^{\gamma},p_{0}^{\gamma}):= n\int_{0}^{\pi}\hspace{-0.25cm}\int  \frac{(p_{m,h}^{\gamma}(z,\phi)-p_{0}^{\gamma}(z,\phi))^2}{p_{0}^{\gamma}(z,\phi)}dzd\phi\leq  \frac{M}{4}.
\end{eqnarray*}
\end{itemize}
\noindent Proofs of this three conditions are done in Appendix~\ref{C1} to \ref{C3}.

\subsubsection{\textbf{Construction of} $\mathbf{W_0}$}
\label{W0}

The Wigner function $W_0$ is the same as in the paper of \cite{Butucea&Guta&Artiles}. For the sake of completeness, we recall its construction here. The probability density function associated to any density matrix $\rho$ in the ideal noiseless setting is given by equation (\ref{de.p.a.rho}).  In particular, for diagonal density matrix $\rho$, the associated probability density function is
\begin{eqnarray*}
         \label{lienprobadensite}
         p_\rho(x,\phi)=\sum_{k=0}^\infty \rho_{kk}\psi_k^2(x).
\end{eqnarray*}
For all $0<\alpha,\lambda<1$, we introduce a family of diagonal density matrix $\rho^{\alpha,\lambda}$ such that for all $k\in\mathbb N$
\begin{eqnarray}
       \label{rhoalphalambda}
       \rho^{\alpha,\lambda}_{kk}\hspace{-0.15cm}&=&\hspace{-0.15cm}\int_0^1 z^k \alpha\frac{(1-z)^\alpha}
       {(1-\lambda)^\alpha}\1_{\lambda\leq z\leq 1}dz.
\end{eqnarray}
Therefore the probability density associated to this diagonal density matrix  $\rho^{\alpha,\lambda}$ can be written as follow
\begin{eqnarray}
       \hspace{-1cm}
       \label{lienprobadensite-2}
       p_{\alpha,\lambda}(x,\phi)=\sum_{k=0}^\infty \rho_{kk}\psi_k^2(x)=\sum_{k=0}^\infty
       \psi_k^2(x) \int_0^1 z^k \alpha\frac{(1-z)^\alpha}{(1-\lambda)^\alpha}\1_{\lambda\leq z\leq
        1}dz.
\end{eqnarray}
Moreover by the well known Mehler formula (see \cite{Erdelyi}), we have
$$
\sum_{k=0}^\infty z^k\psi_k^2(x)=\frac{1}{\sqrt{\pi(1-z^2)}}\exp\left(-x^2\frac{1-z}{1+z}   \right).
$$
Then, it comes
\begin{eqnarray*}\hspace{-1cm}
          p_{\alpha,\lambda}(x,\phi)\hspace{-0.15cm}&=&\hspace{-0.15cm}\frac{\alpha}{(1-\lambda)^\alpha} \int_0^1 \frac{(1-z)^
          \alpha}{\sqrt{\pi(1-z^2)}}\exp\left(-x^2\frac{1-z}{1+z}   \right)   \1_{\lambda\leq z\leq 1}dz.
\end{eqnarray*}
The following Lemma, proved in the paper of \cite{Butucea&Guta&Artiles}, gives a control on the tails of the associated density $p_{\alpha,\lambda}(x,\phi)=p_{\alpha,\lambda}(x)$ as it doesn't depend on $\phi$.

\begin{Lemma}[Butucea, Guta and Artiles (2007)]
\label{p0ltailbound}
For all $\phi\in[0,1]$ and all $0<\alpha,\lambda<1$ and $|x|>1$ there exist constants $c,C$ depending on $\alpha$ and $\lambda$ such that 
$$
c|x|^{-(1+2\alpha)}\leq p_{\alpha,\lambda}(x) \leq C|x|^{-(1+2\alpha)}.
$$
\end{Lemma}

\noindent In view of Lemma 3 of \cite{Butucea&Guta&Artiles}, the Wigner function $W_0$ will be chosen in the set
 $$
 \mathcal{W}^{\alpha,\lambda}=\left\{ W^{\alpha,\lambda}=W_{\rho_{\alpha,\lambda}}  :\rm{ Wigner\, \,function\, \,associated\, \, to \, \, } \rho_{\alpha,\lambda}: 0<\alpha,\lambda<1\right\},
 $$
with $\lambda$ close enough to $1$ so that $W_0 \in \Acal(\beta,L)$ (see \cite{Butucea&Guta&Artiles} for the proof and details).

\subsubsection{\textbf{Construction of the set of Wigner functions $\mathbf{\mathcal W}_{\delta,h}$ for the $\mathbf{\mathbb L_2}$-norm}}
\label{blabla}

\noindent\\
\noindent 
We define $M+1$ infinitely differentiable functions such that: 
\begin{itemize}
\item For all $m=1\cdots,M$, $g_m\,:\, \mathbb R \rightarrow [0,1]$.
\item The support of $g_m$ is 
$
\mathrm{Supp}(g_m) = \left] m\delta,(m+1)\delta \right[.
$
\item And
$
\forall t\in \left[ (m+1/3)\delta,(m+2/3)\delta \right], \, g_m(t)=1.
$
\item An odd function $g\,:\, \mathbb R \rightarrow [-1,1]$, such that for some fixed $\epsilon>0,$ $g(x)=1$ for any $x\geq \epsilon$.
\end{itemize}

\noindent Define also the following constants :  
\begin{eqnarray}
\label{am}
&&a_m := (h^{-2} + m\delta)^{1/2},\quad b_m := (h^{-2} + (m+1)\delta)^{1/2},\quad \forall m=1,\cdots,M.\\
\label{tildeam}
&&\widetilde a_m := (h^{-2} + (m+1/3)\delta)^{1/2}, \quad\widetilde b_m := (h^{-2} + (m+2/3)\delta)^{1/2}\quad \forall m=1,\cdots,M.\\
\label{C0}
&&C_0:=\sqrt{\pi L (\beta+\gamma)}.
\end{eqnarray}
\noindent We also introduce $M$  infinitely differentiable functions such that:

\begin{itemize}
\item For all $m=1\cdots,M$, $V_{m,h}:   \R^2\rightarrow \R$  is an odd real-valued function. 
\item Set $t = \sqrt{w_1^2 + w_2^2}$, then the function $V_{m,h}$ admitting Fourier transform with respect to both variable equals to
\begin{eqnarray}
\label{jh}
\widetilde{V}_{m,h}(w):= \mathcal F_2[V_{m,h}](w):=i a C_0h^{-1}e^{\beta h^{-2}}e^{-2\beta |t|^{2}}g_{m}(|t|^2-h^{-2})g(w_2),
\end{eqnarray}
where $a>0$ is a numerical constant chosen sufficiently small. The bandwidth is such that
\begin{eqnarray}
       \label{h2}
h \hspace{-0.15cm}&=&\hspace{-0.15cm} \left( \frac{\log n}{2(\beta+ \gamma)}\right)^{-1/2}.
\end{eqnarray}
\end{itemize}

\noindent Note that $\widetilde{V}_{m,h}(w)$ is infinitely differentiable and compactly supported, thus it belongs to the Schwartz class $\mathcal S(\R^2)$ of fast decreasing functions on $\R^2$. The Fourier transform being a continuous mapping of the Schwartz class onto itself, this implies that $V_{m,h}$ is also in the Schwartz class $\mathcal S(\R^2)$. Moreover, $\widetilde{V}_{m,h}(w)$ is an odd function with purely imaginary values. Consequently, $V_{m,h}$ is an odd real-valued function. Consequently, we get
\begin{eqnarray}
       \hspace{-1cm}\label{int-Vh}
\iint V_{m,h}(p,q)dpdq=\int \mathcal R[V_{m,h}](x,\phi)dx=0,
\end{eqnarray}
 for all $\phi\in[0,\pi]$ and $\mathcal R[V_{m,h}]$ the Radon transform of $V_{m,h}$.  As in  (\ref{eq:prhophi}),  we define 
\begin{eqnarray}
p_{m,h}(x,\phi)=\frac{1}{\pi}\mathcal R[W_{m,h}] (x,\phi) \1_{(0,\pi}(\phi),\nonumber\\
\label{bienMD}
\text{and}\quad\rho^{(m,h)}_{j,k}=\int_{0}^{\pi}\hspace{-0.25cm}\int p_{m,h}(x,\phi) f_{j,k}(x)e^{(j-k)\phi}dxd\phi.
\end{eqnarray}
By Lemma~\ref{bonne matrice} in Appendix~\ref{proofbonne matrice}, the matrix $\rho^{(m,h)}$ is proved to be a density matrix. Therefore, in view of (\ref{de.p.a.rho}) and (\ref{int-Vh}), the function $W_{m,h}$ is a Wigner function.  Now, we can define our set of Wigner functions
\begin{eqnarray}
\label{claaseL2}
\mathcal{W}_{\delta,h}=\left\{ W_{m,h}:   \R^2\rightarrow \R, \, W_{m,h}(z)=W_0(z)+ V_{m,h}(z),\,m=1,\cdots,M\right\},
\end{eqnarray}
where $W_0$ is the Wigner function associated to the density $p_0$ defined in (\ref{lienprobadensite}).

\subsubsection{\textbf{Condition $\mathbf{(C1)}$}}
\label{C1}

\noindent  By the triangle inequality and for any $1\leq m \leq M$,  we have 
\begin{eqnarray*}\hspace{-1cm}
\|  \widetilde{W}_{m,h}e^{\|\cdot\|^2}\|_2\leq \|  \widetilde{W}_{0}e^{\|\cdot\|^2}\|_2+\|  \widetilde{V}_{m,h}e^{\|\cdot\|^2}\|_2.
\end{eqnarray*}

\noindent The first term in the above sum has be bounded in Lemma 3 of \cite{Butucea&Guta&Artiles} as follow
\begin{eqnarray}
\label{B1}
\|  \widetilde{W}_{0}e^{\|\cdot\|^2}\|_2^2 \leq \pi^2 L.
\end{eqnarray}
\noindent To study the second term in the sum above, we consider the change of variables $w=(t\cos \phi,t\sin\phi)$ and as $g$ is bounded by 1, we get  since (\ref{delta}), (\ref{am}) and (\ref{C0}) that
\begin{eqnarray}
\hspace{-1.5cm}
\|  \widetilde{V}_{m,h}e^{\beta\|\cdot\|^2}\|_2^2
\hspace{-0.15cm}&\leq&\hspace{-0.15cm}\iint  \left[a C_0h^{-1}e^{\beta h^{-2}}\right]^2 e^{-2\beta \|w\|^{2}}g_{m}^2(\|w\|^2-h^{-2})dw\nonumber\\
\hspace{-0.15cm}&\leq&\hspace{-0.15cm}a^2 C_0^2 h^{-2}e^{2\beta h^{-2}} \int_0^\pi\int_{a_m}^{b_m} |t| e^{-2\beta |t|^{2}}dt\nonumber\\
\hspace{-0.15cm}&\leq&\hspace{-0.15cm} \pi a^2 C_0^2 h^{-2} e^{2\beta h^{-2}}e^{-2\beta a_m^{2}}\int_{a_m}^{b_m} t dt
\leq \frac{\pi}{2}  a^2 C_0^2h^{-2} e^{-2\beta m\delta}\left[ b_m^2-a_m^2  \right]\nonumber\\
\label{B2}
\hspace{-0.15cm}&\leq&\hspace{-0.15cm}\frac{\pi}{3} a^2 C_0^2 h^{-2}\delta e^{-2\beta m\delta} \leq\pi^2 L,
\end{eqnarray}
for $a$ small enough.  Combining (\ref{B1}) and (\ref{B2}), it comes $W_{m,h} \in \Acal(\beta,L)$ for any $1\leq m \leq M$.

\subsubsection{\textbf{Condition $\mathbf{(C2)}$}}
\label{C2}

\noindent By applying Plancherel Theorem  and the change of variables $w=(t\cos \phi,t\sin\phi)$, we have since the supports $\mathrm{Supp}(g_k)$ and $\mathrm{Supp}(g_m)$ are disjoints for any $k\neq m$ that
\begin{eqnarray}\hspace{-0.1cm}
\|W_{k,h}- W_{m,h} \|_2^2\hspace{-0.15cm}&=&\hspace{-0.15cm}\|V_{k,h}- V_{m,h} \|_2^2=
\frac{1}{4\pi^2}\int_{0}^{\pi}\hspace{-0.25cm}\int |t| \left|   \widetilde{V}_{k,h}(t,\phi)-\widetilde{V}_{m,h}(t,\phi)      \right|^2dtd\phi\nonumber\\
\label{C21}
\hspace{-0.15cm}&=&\hspace{-0.15cm}\frac{a^2C_0^2}{4\pi^2} h^{-2} e^{2\beta h^{-2}}\int_{0}^{\pi}\hspace{-0.25cm}\int |t| e^{-4\beta t^2} g^2(t\sin\phi)\left[  g_{k}^2(t^2 - h^{-2}) + g_{m}^2(t^2 - h^{-2}) \right] dtd\phi.
\end{eqnarray}
\noindent  Note that for a fixed $\mu \in ]0,\pi/4[$, there exists a numerical constant $c>0$ such that $\sin(\phi)>c$ on $]\mu,\pi-\mu[$. From now, we denote by $\widetilde A_m$ the set  
\begin{eqnarray}
\label{Am}\widetilde A_m:=\left\{w\in\R^2:\, (m+1/3)\delta\leq \|w\|^2\leq (m+2/3)\delta\right\},\quad \forall m=1,\cdots,M.
\end{eqnarray} 
By definition of $g$ and for a large enough $n$, we have for any $(t,\phi) \in (\widetilde  A_k \cup \widetilde  A_m) \times ]\mu,\pi-\mu[$ that $g^2(t\sin(\phi)) = 1$ with  $t^2=\|w\|^2$. Therefore, (\ref{C21}) can be lower bounded as follows
\begin{eqnarray}
\|W_{k,h}-W_{m,h} \|_2^2
\hspace{-0.15cm}&\geq&\hspace{-0.15cm}\frac{a^2C_0^2}{4\pi^2} h^{-2} e^{2\beta h^{-2}}\int_{\mu}^{\pi-\mu}\hspace{-0.25cm}\int_{\widetilde A_k \cup \widetilde A_m}|t| e^{-4\beta t^2} \left[  g_{k}^2(t^2 - h^{-2})+g_{m}^2(t^2 - h^{-2}) \right] dtd\phi\nonumber\\
\label{C22}
\hspace{-0.15cm}&=&\hspace{-0.15cm}\frac{\pi - 2\mu}{4\pi^2} a^2C_0^2h^{-2} e^{2\beta h^{-2}} \int_{\widetilde A_k \cup \widetilde A_m}|t| e^{-4\beta t^2} \left[  g_{k}^2(t^2 - h^{-2})+g_{m}^2(t^2 - h^{-2}) \right] dt.
\end{eqnarray}
\noindent On $\widetilde A_m$ and by construction of the function $g_m$, we have     
$$g_{m}^2(t^2 - h^{-2})=1,\quad 1\leq m\leq M.$$
\noindent Constants defined in (\ref{tildeam}) are such that for $k>m$, we have $\widetilde a_{m} <\widetilde b_{m} <\widetilde a_{k} <\widetilde b_{k} $.  Whence, since $\widetilde A_m$ and $\widetilde A_k$ are disjoint sets for any $k>m$, it results 
\begin{eqnarray}\hspace{-1.5cm}
 I\hspace{-0.15cm}&:=&\hspace{-0.15cm}\int_{\widetilde A_k \cup \widetilde A_m}|t| e^{-4\beta t^2} \left[  g_{k}^2(t^2 - h^{-2})+g_{m}^2(t^2 - h^{-2}) \right] dt\nonumber\\
\hspace{-0.15cm}&\geq&\hspace{-0.15cm} e^{-4\beta \widetilde b_{k} ^{2}} \int_{\widetilde A_k \cup\widetilde  A_m}|t|  \left[  g_{k}^2(t^2 - h^{-2})+g_{m}^2(t^2 - h^{-2}) \right] dt\nonumber\\
\label{C23}
 \hspace{-0.15cm}&\geq&\hspace{-0.15cm}
2e^{-4\beta \widetilde b_{k} ^{2}} \int_{\widetilde a_k}^{\widetilde b_k} t dt\geq  e^{-4\beta \widetilde b_{k} ^{2}}(\widetilde b_k^2 - \widetilde a_k^2)t
 \geq\frac{1}{3} \delta e^{-4\beta \widetilde b_k^2}.%
\end{eqnarray}
\noindent Combining (\ref{C22}) and (\ref{C23}), we get since $C_0^2h^{-2}\delta=\pi L/2$ 
\begin{eqnarray*}\hspace{-1cm}
\|W_{k,h}- W_{m,h} \|_2^2 \hspace{-0.15cm}&\geq&\hspace{-0.15cm}\frac{\pi - 2\mu}{12\pi^2} a^2C_0^2h^{-2} e^{2\beta h^{-2}} \delta e^{-4\beta \widetilde b_k^2}=\frac{\pi - 2\mu}{24\pi} a^2 L e^{2\beta h^{-2}} e^{-4\beta \widetilde b_k^2}\\
\hspace{-0.15cm}&=&\hspace{-0.15cm}\frac{\pi - 2\mu}{24\pi} a^2 L e^{-2\beta h^{-2}} e^{-4\beta(k+2/3)\delta}.
\end{eqnarray*}
Since  $1\leq k\leq M\leq 1/\delta$ and(\ref{h2}), it comes
\begin{eqnarray*}\hspace{-2cm}
\|W_{k,h}- W_{m,h} \|_2^2 \hspace{-0.15cm}&\geq&\hspace{-0.15cm}\frac{\pi - 2\mu}{24\pi} a^2 Ln^{-\frac{\beta}{\beta+\gamma}} e^{- 8\beta }\geq  4c  n^{-\frac{\beta}{\beta+\gamma}}=:4\varphi_n^2,
\end{eqnarray*}
where $c>0$ is a numerical constant.

\subsubsection{\textbf{Condition $\mathbf{(C3)}$}}
\label{C3}

\noindent Denote by $\widetilde C>0$  a constant whose value may change from line to line and recall that $N^{\gamma}$ is the density of the Gaussian distribution with zero mean and variance $2\gamma$. Note that $p_0$ and $N^{\gamma}$ do not depend on $\phi$. Consequently, in the framework of noisy data defined in (\ref{Z}), $p_0^{\gamma}(z,\phi) = p_{0}^{\gamma}(z)\frac{1}{\pi}\1_{(0,\pi)}(\phi)$.

\begin{Lemma}
\label{fact1}
There exists numerical constants $c'>0$ and $c''>0$ such that
\begin{eqnarray}
       \hspace{-1.cm}
       \label{tail-conv1}
p_0^{\gamma}(z) \geq c' z^{-2},\quad \forall |z|\geq 1+ \sqrt{2\gamma},
\end{eqnarray}
and
\begin{eqnarray}
       \hspace{-1cm}\label{tail-conv2}
p_0^{\gamma}(z) \geq c'' ,\quad \forall |z|\leq 1+ \sqrt{2\gamma}.
\end{eqnarray}
\end{Lemma}
\noindent The proof of this Lemma is done in Appendix~\ref{prooffact1}. Using Lemma~\ref{fact1}, the  $\chi^2$-divergence can be upper bounded as follow
\begin{eqnarray}
       \hspace{-1.5cm}
       \label{xi-bound}
n\mathcal X^2(p_{m,h}^{\gamma},p_{0}^{\gamma})\hspace{-0.15cm}&=&\hspace{-0.15cm}n\int_{0}^{\pi}\hspace{-0.25cm}\int  \frac{\left(p_{m,h}^{\gamma}(z,\phi)-p_{0}^{\gamma}(z,\phi)\right)^2}{p_{0}^{\gamma}(z,\phi)}dzd\phi\nonumber\\
\hspace{-0.15cm}&\leq&\hspace{-0.15cm} \frac{n}{c''} \int_{0}^{\pi}\hspace{-0.25cm}\int_{-(1+\sqrt{2\gamma})}^{1+\sqrt{2\gamma}}\left( p_{m,h}^{\gamma}(z,\phi)-p_{0}^{\gamma}(z,\phi)  \right)^2 dzd\phi  \nonumber \\
&& 
+ \frac{n}{c'} \int_{0}^{\pi}\hspace{-0.25cm}\int_{\mathbb R\setminus (1+\sqrt{2\gamma},1+\sqrt{2\gamma})} z^2\left( p_{m,h}^{\gamma}(z,\phi)-p_{0}^{\gamma}(z,\phi)  \right)^2 dzd\phi \nonumber \\
\hspace{-0.15cm}&=:&\hspace{-0.15cm}\frac{n}{c''} I_1+\frac{n}{c'} I_2.
\end{eqnarray}

\noindent First underline, as in (\ref{fourierproun})  the Fourier transforms of $p_{m,h}^{\gamma}$ and $p_0^{\gamma}$ with respect to the first variable are equal respectively to 
\begin{eqnarray}
\label{fourierdepmh}
\Fcal_1[p_{m,h}^{\gamma}(\cdot,\phi)](t)\hspace{-0.15cm}&=&\hspace{-0.15cm}\widetilde W_{m,h}(t\cos\phi,t\sin\phi)\widetilde{N}^\gamma(t)\nonumber\\
\hspace{-0.15cm}&=&\hspace{-0.15cm}\left(\widetilde V_{m,h}(t\cos\phi,t\sin\phi)+ \widetilde W_0(t\cos\phi,t\sin\phi)\right)e^{-\gamma t^2},\\
\label{fourierdep0}
\Fcal_1[p_{0}^{\gamma}(\cdot,\phi)](t)\hspace{-0.15cm}&=&\hspace{-0.15cm}\widetilde W_0(t\cos\phi,t\sin\phi) e^{-\gamma t^2},
\end{eqnarray}
\noindent since  $\widetilde{N}^\gamma(t)=e^{-\gamma t^2}$. Using Plancherel Theorem and (\ref{jh}), equations (\ref{fourierdepmh}) and (\ref{fourierdep0}), the first integral $I_1$ in the sum (\ref{xi-bound}) is bounded by
\begin{eqnarray*}
       \hspace{-1cm}
I_1\hspace{-0.15cm}&\leq&\hspace{-0.15cm}\int_{0}^{\pi}\hspace{-0.25cm}\int \left( p_{m,h}^{\gamma}(z,\phi)-p_{0}^{\gamma}(z,\phi)  \right)^2 dzd\phi 
=\frac{1}{4\pi^2}\int_{0}^{\pi}\hspace{-0.25cm}\int  \left| \mathcal F_1[ p_{m,h}^{\gamma}(\cdot,\phi)](t)-\mathcal F_1[p_{0}^{\gamma}(\cdot,\phi)] (t) \right|^2 dtd\phi \nonumber\\
\hspace{-0.15cm}&=&\hspace{-0.15cm} \frac{1}{4\pi^2}\int_{0}^{\pi}\hspace{-0.25cm}\int  \left|\widetilde{V}_{m,h}(t\cos\phi,t\sin\phi)\right|^2e^{-2\gamma t^2} dtd\phi \nonumber\\
\hspace{-0.15cm}&=&\hspace{-0.15cm}\frac{a^2 C_0^2}{4\pi^2} h^{-2} e^{2\beta h^2}\int_{0}^{\pi}\hspace{-0.25cm}\int e^{-4\beta t^2-2 \gamma t^2} g_m^2\left(t^2 - h^{-2} \right)g^2(t\sin \phi) dtd\phi.\nonumber
\end{eqnarray*}
By construction, the function $g$ is bounded by 1 and the function $g_m$ admits as support $\text{Supp}(g_m)=]m\delta,(m+1)\delta[$ for all $m=1,\cdots,M$. Thus, 
\begin{eqnarray*}
       \hspace{-1cm}
I_1
\hspace{-0.15cm}&\leq&\hspace{-0.15cm}\frac{a^2 C_0^2}{4\pi} e^{2\beta h^2} \int e^{-4\beta t^2- 2 \gamma t^2} g_m^2\left(t^2 - h^{-2} \right) dt
\leq \frac{a^2 C_0^2}{4\pi} h^{-2} e^{2\beta h^{-2}} \int_{a_m}^{b_m} e^{-4\beta t^2- 2 \gamma t^2}dt\nonumber\\
\hspace{-0.15cm}&\leq&\hspace{-0.15cm} \frac{a^2 C_0^2}{4\pi} (b_m - a_m)h^{-2} e^{2\beta h^{-2}} e^{-4\beta a_m^2- 2 \gamma a_m^2}
\leq \frac{a^2 C_0^2}{4\pi} \frac{b_m^2 - a_m^2}{2 a_m} h^{-2} e^{2\beta h^{-2} -4\beta a_m^2- 2 \gamma a_m^2}.\nonumber
\end{eqnarray*}
 Some basic algebra, (\ref{delta}), (\ref{am}), (\ref{C0}) and (\ref{h2}) yield
\begin{eqnarray}
       \hspace{-1cm}\label{xi-bound2}
\frac{n}{c''} I_1 \leq \frac{a^2 \widetilde  C}{\sqrt{\log n}},
\end{eqnarray}
for some a constant $\widetilde  C>0$ whose may depend on $\beta$, $\gamma$, $L$ and $c''$. For the second term $I_2$ in the sum (\ref{xi-bound}), with the same tools we obtain using  in addition the spectral representation of the differential operator, that 
\begin{eqnarray}\hspace{-1cm}
\label{I2}
I_2\hspace{-0.15cm}&\leq&\hspace{-0.15cm}\int_{0}^{\pi}\hspace{-0.25cm}\int z^2 \left( p_{m,h}^{\gamma}(z,\phi)-p_{0}^{\gamma}(z,\phi)  \right)^2 dzd\phi \nonumber \\
\hspace{-0.15cm}&=&\hspace{-0.15cm}\int_{0}^{\pi}\hspace{-0.25cm}\int \left|\frac{\partial}{\partial t}\left(  \Fcal_1[ p_{m,h}^{\gamma}(\cdot,\phi)]-\mathcal F_1[p_{0}^{\gamma}(\cdot,\phi)]  \right)(t)\right|^2 dt d\phi\nonumber \\
\hspace{-0.15cm}&=&\hspace{-0.15cm}\int_{0}^{\pi}\hspace{-0.25cm}\int  \left|\frac{\partial}{\partial t}\left(  \widetilde{V}_{m,h}(t\cos\phi,t\sin\phi) e^{-\gamma t^2} \right)\right|^2 dtd\phi \nonumber  \\
\hspace{-0.15cm}&=&\hspace{-0.15cm}\int_{0}^{\pi}\hspace{-0.25cm}\int \left| e^{-\gamma t^2}\frac{\partial}{\partial t}(\widetilde{V}_{m,h})(t\cos \phi,t \sin \phi)  -2\gamma te^{-\gamma t^2}  \widetilde{V}_{m,h}(t\cos \phi,t\sin \phi) \right|^2 dt d\phi\nonumber\\
\hspace{-0.15cm}&\leq&\hspace{-0.15cm}2\int_{0}^{\pi}\hspace{-0.25cm}\int e^{-2\gamma t^2} \left| I_{2,1} \ \right|^2 dtd\phi+16\gamma^2\int_{0}^{\pi}\hspace{-0.25cm}\int t^2e^{-2\gamma t^2}    \left| I_{2,1}\right|^2 dtd\phi,
\end{eqnarray}
where $I_{2,2}=\widetilde{V}_{m,h}(t\cos \phi,t\sin \phi)$ and $I_{2,1}$, the partial derivative $\frac{\partial}{\partial t}(\widetilde{V}_{m,h})(t\cos \phi,t\sin \phi)$, is equal to
$$
ia C_0 h^{-1}e^{\beta h^{-2}-2\beta t^2} \left[ g_m(t^2 - h^{-2})\left( -4\beta tg(t\sin\phi) +g'(t\sin\phi) \sin \phi\ \right)+ 2t g_m'(t^2 - h^{-2})g(t\sin \phi) \right].
$$
Since $g_m$ and $g$ belong to the Schwartz class, there exists a numerical constant $c_S>0$ such that $\max\{\|g_m\|_{\infty}, \|g_m'\|_{\infty}, \|g\|_{\infty}, \|g'\|_{\infty}\} \leq c_S$.  Furthermore, for all $m=1,\cdots,M$, the support of the function $g_m$ is  $\text{Supp}(g_m)=]m\delta,(m+1)\delta[$, then
\begin{eqnarray}
\label{I21}
 \left| I_{2,1}  \right|^2\hspace{-0.15cm}&\leq&\hspace{-0.15cm}
a^2 c_S^4C_0^2 h^{-2}e^{2\beta h^{-2}-4\beta t^2}   \left((4\beta +2)|t| +1 \right)^2\1_{(a_m,b_m)}(t),%
\end{eqnarray}
\noindent with $a_m$ and $b_m$ defined in (\ref{am}). Similary, we have
\begin{eqnarray}
\label{I22}
\left| I_{2,2} \right|^2\hspace{-0.15cm}&=&\hspace{-0.15cm}  \left| a C_0 h^{-1}e^{\beta h^{-2}} e^{-2\beta t^2} g_m(t^2 - h^{-2})g(t\sin\phi)\right|^2\nonumber\\
\hspace{-0.15cm}&\leq&\hspace{-0.15cm}
a^2 c_S^4C_0^2 h^{-2}e^{2\beta h^{-2}-4\beta t^2}   \1_{(a_m,b_m)}(t).\nonumber\\
\end{eqnarray}
\noindent Combining (\ref{I21}) and (\ref{I22}) with (\ref{I2}), as $0\leq m\delta\leq 1$
\begin{eqnarray*}\hspace{-1cm}
I_2\hspace{-0.15cm}&\leq&\hspace{-0.15cm}2a^2 c_S^4C_0^2 h^{-2}e^{2\beta h^{-2}}
\int_{0}^{\pi}\hspace{-0.25cm}\int_{a_m}^{b_m} e^{-2\gamma t^2} e^{-4\beta t^2}   \left[\left((4\beta +2)|t| +1 \right)^2+ 8\gamma^2 t^2\right]
dtd\phi\nonumber\\
\hspace{-0.15cm}&\leq&\hspace{-0.15cm}2\pi a^2 c_S^4C_0^2 h^{-2}e^{2\beta h^{-2}}e^{-(4\beta+2\gamma) a_m^2}\left[\left((4\beta +2)b_m +1 \right)^2+ 8\gamma^2 b_m^2\right]\int_{a_m}^{b_m}  dt\nonumber\\
\hspace{-0.15cm}&\leq&\hspace{-0.15cm}2\pi a^2c_S^4C_0^2 h^{-2}e^{-2(\beta + \gamma) h^{-2}}   e^{-2(2\beta+\gamma)m\delta}\left[\left((4\beta +2)b_m +1 \right)^2+ 8\gamma^2 b_m^2\right]\frac{b_m^2 - a_m^2}{2 a_m}\nonumber\\
\hspace{-0.15cm}&\leq&\hspace{-0.15cm}2\pi a^2 c_S^4C_0^2 h^{-2}e^{-2(\beta + \gamma) h^{-2}}   \left[\left((4\beta +2)b_m +1 \right)^2+ 8\gamma^2 b_m^2\right]\frac{\delta}{2 a_m}.\nonumber\\
\end{eqnarray*}

\noindent Some basic algebra, (\ref{delta}), (\ref{am}), (\ref{C0}) and (\ref{h2}) yield
\begin{eqnarray}
\label{I2bis}
\frac{n}{c'}I_2\hspace{-0.15cm}&\leq&\hspace{-0.15cm}  a^2 \widetilde  C\sqrt{\log n},
\end{eqnarray}
for some a constant $\widetilde  C>0$ whose may depend on $\beta$, $\gamma$, $L$ $c_S$ and $c'$.  Combining (\ref{I2bis}) and  (\ref{xi-bound2}) with (\ref{xi-bound}), we get for $n$ large enough
\begin{eqnarray*}\hspace{-1cm}
n\mathcal X^2(p_{k,h}^{\gamma},p_{0}^{\gamma}):= n\int_{0}^{\pi}\hspace{-0.25cm}\int _{\mathbb R} \frac{(p_{k,h}^{\gamma}(z,\phi)-p_{0}^{\gamma}(z,\phi))^2}{p_{0}^{\gamma}(z,\phi)}dzd\phi\leq  a^2 \widetilde  C\sqrt{\log n},
\end{eqnarray*}
where $\widetilde  C>0$ is a constant whose may depend on $\beta$, $\gamma$, $L$ $c_S$, $c"$ and $c'$. Taking the numerical constant $a>0$ small enough, we deduce from the previous display that
$$
n\mathcal X^2(p_{k,h}^{\gamma},p_{0}^{\gamma})\leq \frac{M}{4},
$$
since $M = \lfloor \sqrt{\log n} \rfloor$.

\subsection{Proof of Theorem~2 - Lower bounds for the sup-norm}
\label{proofLBU}

\noindent To prove the lower bound for the sup-norm, we need to slightly modify the construction of the Wigner classe $\mathcal{W}_{\delta,h}$ defined in (\ref{claaseL2}) into
\begin{eqnarray}
\label{claaseL2bis}
\mathcal{W}_{\delta,h,\epsilon}=\left\{ W_{m,h,\epsilon}:   \R^2\rightarrow \R, \, W_{m,h,\epsilon}(z)=W_0(z)+ V_{m,h,\epsilon}(z),\,m=1,\cdots,M\right\},
\end{eqnarray}
where $W_0$ is the Wigner function associated to the density $p_0$ defined in (\ref{lienprobadensite}) stay unchanged as compared to the $L_2$ case. However, the construction of the $\{V_{m,h}\}_{m}$-functions defined in (\ref{jh}) only changed through modification of the functions $g_m$ and $g$ respectively into $g_{m,\epsilon}$ and $g_{\epsilon}$, for $0<\epsilon<1$.

\noindent\\
\noindent 
We define $M+1$ infinitely differentiable functions such that: 
\begin{itemize}
\item For all $m=1\cdots,M$, $g_{m,\epsilon}\,:\, \mathbb R \rightarrow [0,1]$.
\item The support of $g_{m,\epsilon}$ is 
$
\mathrm{Supp}(g_{m,\epsilon}) = \left] m\delta,(m+1)\delta \right[.
$
\item Using a similar construction as for function $g_m$, we can also assume that
\begin{eqnarray}
       \hspace{-1cm}\label{der-g}
g_{m,\epsilon}(t)  = 1, \quad \forall t\in B_{m,\epsilon} := \left[( m+\epsilon)\delta, ( m+1-\epsilon)\delta \right],
\end{eqnarray}
\noindent and
\begin{eqnarray}
       \hspace{-1cm}\label{der-gm}
\|g_{m,\epsilon}'\|_{\infty} \leq \frac{c}{\epsilon \delta},
\end{eqnarray}
for some numerical constant $c>0$.

\item An odd function $g_{\epsilon}\,:\, \mathbb R \rightarrow [-1,1]$ satisfies the same conditions as $g$ above but we assume in addition that
\begin{eqnarray}
       \hspace{-1cm}\label{geps}
\|g_{\epsilon}'\|_{\infty} \leq \frac{c}{\epsilon},
\end{eqnarray}
for some numerical constant $c>0$.
\end{itemize}
\noindent The condition (\ref{geps}) will be needed to check Condition \textbf{(C3)}. Such a function can be easily constructed. Consider for instance a function $g_{\epsilon}$ such that its derivative satisfies 
$$g_{\epsilon}'(t) =  \left[\psi\ast\frac{1}{\epsilon}\1_{(0,\epsilon)}\right](t),$$
\noindent  for any $t\in (0,\epsilon)$ where $\psi$ is a mollifier. Integrate this function and renormalize it properly so that $g_{\epsilon}(t) =1$ for any $t\geq \epsilon$. Complete the function by symmetry to obtain an odd function defined on the whole real line. Such a construction satisfies condition (\ref{geps}).\\

\noindent  It is easy to see that Condition \textbf{(C1)} is always satisfied by the new test functions $\{W_{m,h,\epsilon}\}_m$. To check Condition \textbf{(C2)} set $C_h = i a C_0 h^{-1}e^{\beta h^{-2}}$ and then we have
\begin{eqnarray*}\hspace{-1cm}
W_{k,h,\epsilon}(z) -W_{m,h,\epsilon}(z) \hspace{-0.15cm}&=&\hspace{-0.15cm}\frac{1}{4\pi^2} \iint e^{- i\langle  z, w \rangle} \left(\widetilde W_{k,h,\epsilon}(w) - \widetilde W_{m,h,\epsilon}(w)\right)dw\\
\hspace{-0.15cm}&=&\hspace{-0.15cm}\frac{1}{4\pi^2} \int_{0}^{\pi}\hspace{-0.25cm}\int e^{-it[z,\phi]}|t|\left(\widetilde W_{k,h,\epsilon}(t\cos\phi,t\sin\phi)\right.\\
&&\left. - \widetilde W_{m,h,\epsilon}(t\cos\phi,t\sin\phi)\right)dtd\phi\\
\hspace{-0.15cm}&=&\hspace{-0.15cm}\frac{1}{4\pi^2}\int_{0}^{\pi}\hspace{-0.25cm}\int e^{-it[z,\phi]}|t|C_h e^{-2 \beta t^2} (g_{k,\epsilon}-g_{m,\epsilon})\left( t^2 - h^{-2} \right) g_{\epsilon}(t)dtd\phi.
\end{eqnarray*}
For all $z\in\R^2$ and $B_{m,0}=\lim_{\epsilon \rightarrow 0}B_{m,\epsilon}$ defined in (\ref{der-g}), we define the following quantity
$$
I(z): = \int_{0}^{\pi}\hspace{-0.25cm}\int e^{-it[z,\phi]}|t|C_h e^{-2 \beta t^2} [\1_{B_{k,0}}-\1_{B_{m,0}}]\left( t^2 - h^{-2} \right)\left[\1_{(0,\infty)}(t) - \1_{(-\infty,0)}(t) \right] dtd\phi.
$$
Lebesgue dominated convergence Theorem guarantees that 
\begin{eqnarray*}\hspace{-1cm}
&&\lim_{\epsilon \rightarrow 0} \left(\int_{0}^{\pi}\hspace{-0.25cm}\int e^{-it[z,\phi]}|t|C_h e^{-2 \beta t^2} (g_{k,\epsilon}-g_{m,\epsilon})\left( t^2 - h^{-2}  \right)g_{\epsilon}(t)dtd\phi \right)=I(z).
\end{eqnarray*}
Therefore, there exists an $\epsilon>0$ (possibly depending on $n,z$) such that
$$
\left|W_{k,h,\epsilon}(z) -W_{m,h,\epsilon}(z)\right| \geq \frac{1}{2} \left| I(z)  \right|.
$$
Taking $z=(0,2 h)$, Fubini's Theorem gives
\begin{eqnarray*}\hspace{-0.5cm}
I(z) \hspace{-0.15cm}&=&\hspace{-0.15cm}\frac{1}{4\pi^2} \int_{0}^{\pi}\hspace{-0.25cm}\int e^{-it2h\sin\phi}|t|C_h e^{-2 \beta t^2} [\1_{A_{k,0}}-\1_{A_{m,0}}]\left( t^2 - h^{-2}   \right)\\
&&\hspace{0.15cm}\times\left[\1_{(0,\infty)}(t) - \1_{(-\infty,0)}(t) \right]dtd\phi\\
\hspace{-0.15cm}&=&\hspace{-0.15cm} \frac{1}{4\pi^2}  \int\hspace{-0.15cm} \left( \int_{0}^{\pi}e^{-it2h\sin\phi}d\phi\right)|t|C_h e^{-2 \beta t^2} [\1_{A_{k,0}}-\1_{A_{m,0}}]\left( t^2 - h^{-2}   \right)\\
&&\hspace{0.15cm}\times\left[\1_{(0,\infty)}(t) - \1_{(-\infty,0)}(t) \right]dt.
\end{eqnarray*}
\noindent Note that
$$
\int_{0}^{\pi} e^{-it2h\sin\phi}d\phi = \pi (iH_{0}(2ht) + J_{0}(2ht)),
$$
where $H_0$ and $J_0$ denote  respectively the Struve and Bessel functions of order $0$. By definition, $H_0$ is an odd function while $J_0$ and $t\rightarrow |t|C_h e^{-2 \beta t^2} [\1_{A_{k,0}}-\1_{A_{m,0}}]\left(t^2 - h^{-2} \right)$ are even functions. Consequently, we get
\begin{eqnarray*}\hspace{-0.5cm}
I(z) &=  &\frac{1}{4\pi} i C_h \int  |t| H_{0}(2ht) e^{-2 \beta t^2} [\1_{A_{k,0}}-\1_{A_{m,0}}]\left( t^2 - h^{-2}   \right)\left[\1_{(0,\infty)}(2ht) - \1_{(-\infty,0)}(t) \right]dt\\
\hspace{-0.15cm}&=&\hspace{-0.15cm}\frac{1}{2\pi} i C_h \int_{0}^{\infty} t H_{0}(2ht) e^{-2 \beta t^2} [\1_{A_{k,0}}-\1_{A_{m,0}}]\left( t^2 - h^{-2}   \right)dt\\
 \hspace{-0.15cm}&=&\hspace{-0.15cm}\frac{iC_h }{2\pi} \left(\int_{a_k}^{b_k} t H_0(2h t) e^{-2 \beta t^2}dt - \int_{a_m}^{b_m} tH_0(2h t) e^{-2 \beta t^2} dt]\right),
\end{eqnarray*}
with $a_k$ and $b_k$ defined in (\ref{am}). For some numerical constant $c>0$, 
$$
]a_m,b_m[ \subset [h^{-1}, h^{-1} + ch \delta].
$$
\noindent On $[h^{-1}, h^{-1} + ch \delta]$ and for a large enough $n$, functions $t\rightarrow H_0(2 h t )$ and $t\rightarrow t e^{-2\beta t^2}$ are decreasing and  
$$\min_{t\in [h^{-1}, h^{-1} + ch \delta]}\{H_0(2 h t )\} \geq 1/2.$$
\noindent  Assume without loss of generality that $k<m$. We easily deduce from the previous observations that
\begin{eqnarray*}\hspace{-1cm}
|I(z)| \hspace{-0.15cm}&\geq&\hspace{-0.15cm} \frac{|C_h| }{4\pi} \left(\int_{a_k}^{b_k} t e^{-2 \beta t^2}dt - \int_{a_m}^{b_m} t e^{-2 \beta t^2} dt\right)\\
\hspace{-0.15cm}&\geq&\hspace{-0.15cm}\frac{|C_h| }{16 \pi\beta} \left(e^{-2 \beta a_k^2} - e^{-2 \beta b_k^2} + e^{-2 \beta b_m^2} - e^{-2 \beta a_m^2}\right)\\
\hspace{-0.15cm}&\geq&\hspace{-0.15cm} \frac{|C_h| }{16 \pi\beta} e^{-2 \beta h^{-2}} e^{-2\beta \alpha m \delta} (1- e^{-2 \beta \alpha (m-k)\delta} )(1 - e^{-2 \beta \delta}).
\end{eqnarray*}
\noindent Therefore, some simple algebra gives 
$$
|I(z)| \geq c \delta^2 |C_h| n^{-\frac{\beta}{\beta+\gamma}} \geq a c' n^{-\frac{\beta}{2(\beta+\gamma)}} \log^{-3/2}(n),
$$
for some numerical constants $c,c'>0$ depending only $\beta$. Taking the numerical constant $a>0$ small enough independently of $n,\beta,\gamma$, we get that Condition \textbf{(C2)} is satisfied with $\varphi_n = c n^{-\frac{\beta}{2(\beta+\gamma)}} \log^{-3/2}(n)$.
 
\noindent Concerning Condition \textbf{(C3)}, we proceed similarly as above for the quadratic risk. The only modification appears in (\ref{I21})-(\ref{I22}) where we now use (\ref{der-g})-(\ref{der-gm}) combined with the fact that 
\begin{center}$|\mathrm{Supp}(g_\epsilon')| \leq 2\epsilon$ and $|\mathrm{Supp}(g_{m,\epsilon}')| \leq 2\delta\epsilon$
\end{center}
\noindent by construction of these functions. Therefore, the details will be omitted here.

\section{Proof of Theorem~3 - Adaptation}%
\label{proofthm-adaptation}

\noindent The following Lemma is needed to prove the Theorem \ref{thm-adaptation}.
%
\begin{Lemma}\label{Lepski}
For $\kappa>0$, a constant, let $\mathcal E_{\kappa} $ be the event defined such that
\begin{eqnarray}
       \label{Eevent}
       \mathcal E_{\kappa} = \bigcap_{m=1}^{M} \left\{ \|\widehat{W}^\gamma_{h_m} - \E
       [\widehat{W}^\gamma_{h_m}]\|_{\infty}\leq \kappa e^{\gamma h_{m}^{-2}}r_n(x+\log M)
       \right\}.
\end{eqnarray}
\noindent Therefore, on the event $\mathcal E_{\kappa}$
$$
\|\widehat{W}^\gamma_{h_{\mathbf{\widehat m }}}  - W_{\rho}\|_{\infty}\leq C\min_{1\leq m \leq M}\left\{  h_m^{r/2-1} e^{-\beta h_m^{-r}} +  e^{\gamma h_{m}^{-2}}r_n(x+\log M)\right\},
$$
where $C>0$ is a constant depending only on $\gamma, \beta, L,r,\kappa$ and $\widehat{W}^\gamma_{h_{\mathbf{\widehat m }}} $ is the adaptive estimator with the bandwidth $h_{\mathbf{\widehat m }}$ defined in (\ref{Lep-proc2}). 
\end{Lemma}
\noindent The proof of the previous Lemma is done in~\ref{proofLepski}. For any fixed $m \in \{1,\cdots,M\}$, we have in view of Proposition \ref{Prop-1-sto-dev} that
$$
\mathbb P \left(\|\widehat{W}^\gamma_{h_m} -\E[\widehat{W}^\gamma_{h_m}] \|_{\infty}\leq C e^{\gamma h_m^{-2}}r_n(x) \right) \geq 1- e^{-x},
$$
\noindent where $r_n(x) = \max\left( \sqrt{\frac{1+x}{n}} , \frac{1+x}{n}\right)$. By a simple union bound, we get 
$$
\mathbb P \left(\bigcap_{1\leq m \leq M}\left\{\|\widehat{W}^\gamma_{h_m} -\E[\widehat{W}^\gamma_{h_m}] \|_{\infty} \leq C_2 e^{\gamma {h_m}^{-2}}r_n(x)\right\}\right) \geq 1- Me^{-x}.
$$
Replacing $x$ by $(x+\log M)$, implies
$$
\mathbb P\left(\bigcap_{1\leq m \leq M}\left\{\|\widehat{W}^\gamma_{h_m} -\E[\widehat{W}^\gamma_{h_m}] \|_{\infty} \leq C_2 e^{\gamma {h_m}^{-2}}r_n(x+\log M)\right\}\right) \geq 1- e^{-x}.
$$
For $\kappa>C_2$, we immediately get that $\mathbb P(\mathcal E_\kappa) \geq 1- e^{-x}$ and the result in probability (\ref{p3}) follows by Lemma \ref{Lepski}. To prove the result in expectation (\ref{p4}),  we use the property $\E[Z] = \int_{0}^{\infty}\mathbb P(Z\geq t)dt$, where $Z$ is any positive random variable. We have indeed for any $1\leq m \leq M$ that
\begin{eqnarray*}\hspace{-1cm}
          \Pro \left( \|\widehat{W}^\gamma_{h_{\widehat l}}  - W_{\rho}\|_{\infty} \geq C\left
          ( h_m^{r/2-1} e^{-\frac{\beta}{h_m^{r}}} +  e^{\gamma h_{m}^{-2}}r_n(x+\log M)\right) 
          \right)\leq e^{-x},\quad \forall x>0.
\end{eqnarray*}
Note that
\begin{eqnarray*}\hspace{-1cm}
                r_n(x+\log M) \hspace{-0.15cm}&=&\hspace{-0.15cm} \max\left\{ \sqrt{\frac{x+\log(eM)}{n}},\frac{x+\log(eM)}{n}\right\}\\
                \hspace{-0.15cm}&\leq&\hspace{-0.15cm} \max\left\{ \sqrt{\frac{\log eM}{n}}, \frac{\log eM}{n}\right\} + \max\left\{\sqrt{\frac
                {x}{n}}\vee\frac{x}{n} \right\}\\
                \hspace{-0.15cm}&\leq&\hspace{-0.15cm} r_n(\log M) + r_n(x-1).
\end{eqnarray*}
\noindent Combining the two previous displays, we get $ \forall x>0$
\begin{eqnarray*}\hspace{-1.5cm}
            \Pro \left( \|\widehat{W}^\gamma_{h_{\widehat l}}  - W_{\rho}\|_{\infty} \geq C\left
            ( h_m^{r/2-1} e^{-\frac{\beta}{h_m^{r}}} +  e^{\gamma h_{m}^{-2}}\left[ r_n(\log M) + r_n
            (x-1)\right]\right) \right)\leq e^{-x}.
\end{eqnarray*}

\noindent Set $Y = \|\widehat{W}^\gamma_{h_{\widehat l}}  - W_{\rho}\|_{\infty}/C$, $a= h_m^{r/2-1} e^{-\frac{\beta}{h_m^{r}}} + e^{\gamma h_{m}^{-2}} r_n(\log M)$ and $b=e^{\gamma h_{m}^{-2}}$. We have
$$
\E[Y] = a+ \E[Y-a] = a+ \int_{0}^{\infty} \mathbb P \left( Y-a \geq u \right) du = a+ b\int_{0}^{\infty} \mathbb P \left( Y-a \geq bt \right) dt.
$$
Set now $t = r_n(x-1)$. If $0< t < 1$, then we have $t = \sqrt{\frac{x}{n}}$. If $t\geq 1$ then we have $t = \frac{x}{n}$. Thus we get by the change of variable $t = \sqrt{\frac{x}{n}}$ that
$$
\int_{0}^{1} \mathbb P \left( Y-a \geq bt \right)dt = \int_{0}^{n}\mathbb P \left( Y-a \geq b\sqrt{\frac{x}{n}} \right)\frac{1}{2\sqrt{x n}}dx \leq  \frac{1}{2\sqrt{n}}\int_{0}^{n}\frac{e^{-x}}{\sqrt{x}}dx \leq \frac{c}{\sqrt{n}},
$$
where $c>0$ is a numerical constant.
Similarly, we get by change of variable $t= \frac{x}{n}$
$$
\int_{1}^{\infty} \mathbb P \left( Y-a \geq bt \right)dt = \int_{n}^{\infty}\mathbb P \left( Y-a \geq b\frac{x}{n} \right)\frac{1}{ n}dx \leq  \frac{1}{n}\int_{n}^{\infty}e^{-x}dx \leq \frac{c'}{n},
$$
where $c'>0$ is a numerical constant. Combining the last three displays, we obtain the result in expectation.

\section{Proof of Auxiliary Lemmas }
\label{sec.proofs.Lemma}

\subsection{\textbf{Proof of Lemma~\ref{VC-class}}}
\label{proofVC-class}
%
\noindent To prove the uniform bound of (\ref{un}), we define
$$
\delta_{h} =\max_{|t|\leq h^{-1}}\left\{|t|e^{\gamma t^2}\right\}.
$$
\noindent  Then, by definition of $K_h^{\eta}$  and by using the inverse Fourier transform formula, we have
\begin{eqnarray}
       \hspace{-1cm}
       \delta_{h}^{-1}\|K_h^{\gamma}\|_{\infty} \hspace{-0.15cm}&=&\hspace{-0.15cm}\frac{1}{2\pi}\delta_{h}^{-1} \sup_{x\in \R} \left|\int e^{-itx} \Fourier{K}_h^{\gamma}(t)dt  \right|  
       \leq \frac{1}{2\pi}\delta_{h}^{-1} \int_{h^{-1}}^{h^{-1}} |t|e^{\gamma t^2}  dt \nonumber \\
       \hspace{-0.15cm}&\leq&\hspace{-0.15cm} \frac{1}{\pi}\delta_{h}^{-1} \int_{0}^{h^{-1}} te^{\gamma t^2}  dt 
       \leq \frac{1}{2 \gamma \pi}\delta_{h}^{-1} \int_{0}^{h^{-1}} 2\gamma te^{\gamma t^2}  dt \nonumber\\
       \label{U}
       \hspace{-0.15cm}&\leq&\hspace{-0.15cm} \frac{1}{2 \gamma \pi}\delta_{h}^{-1} ( e^{\gamma h^{-2}}-1)
       \leq \frac{1}{2 \gamma \pi}\delta_{h}^{-1} ( e^{\gamma h^{-2}}-1)  \leq \frac{h}{2 \gamma \pi}:=U.
\end{eqnarray}

\noindent  For the entropy bound (\ref{deux}), we need to prove that $ K_h^{\ \label{U}}(\cdot)$ admits finite quadratic variation, \textit{i.e.}  $K_h^{\gamma}\in V_2(\R)$, where $V_2(\R)$ is the set of functions with finite quadratic variation (see Theorem 5 of \cite{Bourdaudal2006}). To do this, it is enough to  verify that $K_h^{\gamma} \in B_{2,1}^{1/2}(\R)$ and the result is a consequence  of the embedding $B_{2,1}^{1/2}(\R) \subset V_2(\R)$.\\

\noindent Let us define the Littlewood-Paley characterization of the seminorm $\|\cdot\overset{\bullet}{\|}_{1/2,2,1}$ as follow
\begin{eqnarray*}%
      \| g\overset{\bullet}{\|}_{1/2,2,1} &:=&\sum_{l \in \mathbb Z} 2^{l/2} \|\Fcal_1^{-1}[\alpha_l \Fcal_1[g]]\|_2,
\end{eqnarray*}
\noindent where $\alpha_l(\cdot)$ is a dyadic partition of unity with $\alpha_l$ symmetric w.r.t to $0$, supported in 
$$[-2^{l+1},-2^{l-1}]\cup [2^{l-1}, 2^{l+1}]$$
\noindent  and $0\leq \alpha_1 \leq 1$ (see e.g. Theorem 6.3.1 and Lemma 6.1.7 in the paper of \cite{BerghLof76}). Then, $K_h^{\gamma} \in B_{2,1}^{1/2}(\R)$, if and only if 
$\|K_{h}^{\gamma}\overset{\bullet}{\|}_{1/2,2,1}$  is bounded by a fixed constant. By isometry of the Fourier transform combining with definition of $\alpha_l$ and $K_{h}^{\gamma}$, we get that
\begin{eqnarray*}%
        \|\Fcal_1^{-1}[\alpha_l \Fcal_1[K_{h}^{\gamma}]]\|_{2} \hspace{-0.15cm}&=&\hspace{-0.15cm} \| \alpha_l \Fcal_1[K_{h}^
        {\gamma}] \|_2= \| \alpha_l\Fourier{K}_{h}^{\gamma} \|_2\\\
        \hspace{-0.15cm}&=&\hspace{-0.15cm}\sqrt{2\int_{[0,h^{-1}]\cap [2^{l-1},2^{l+1}]} \alpha_l(t)^2 |t|^2e^{2\gamma t^2}dt}\\
        \hspace{-0.15cm}&\leq&\hspace{-0.15cm}\sqrt{2\int_{[0,h^{-1}]\cap [2^{l-1},2^{l+1}]} t^2 e^{2\gamma t^2}dt}.
\end{eqnarray*}
\noindent A primitive of $t\rightarrow t^{2}e^{2\gamma t^2}$ is $\frac{1}{2\gamma} t e^{2\gamma t^2} - \frac{1}{2\gamma} \int_0^{t} e^{2\gamma u^2}du$. Thus, we get that
\begin{eqnarray*}%
          \|\Fcal_1^{-1}[\alpha_l \Fcal_1[K_{h}^{\gamma}]]\|_{2} &\leq \sqrt{\frac{1}{\gamma}}h^
          {-1/2}e^{\gamma h^{-2}},\quad \forall l \in \mathbb Z,
\end{eqnarray*}
\begin{eqnarray*}%
  \text{and  }\quad   \|K_{h}^{\gamma}\overset{\bullet}{\|}_{1/2,2,1}\leq \sqrt{\frac{1}{\gamma}} h^{-1/2}e^{\gamma h^{-2}}
      \sum_{l =-\infty}^{L_h} 2^{l/2},
\end{eqnarray*}
\noindent where $L_h = \lfloor  \log_2(h^{-1}) +1  \rfloor$.  A simple computation gives that
$$
\sum_{l =-\infty}^{L_{h}} 2^{l/2} \leq \frac{\sqrt{2}}{\sqrt{2}-1}+ \frac{2^{(L_h+1)/2} -1}{\sqrt{2}-1} \leq \frac{\sqrt{2}}{\sqrt{2}-1} + \frac{2}{\sqrt{2}-1}h^{-1/2}.
$$
\noindent  Combining the last two displays and since $h^{-1}\geq 1$, we get
\begin{eqnarray*}
       \|K_{h}^{\gamma}\overset{\bullet}{\|}_{1/2,2,1}\leq c \sqrt{\frac{1}{\gamma}} h^{-1}e^{\gamma h^{-2}},
 \end{eqnarray*}
\noindent where $c>0$ is a numerical constant. This shows that $\delta_h^{-1}\|K_{h}^{\gamma}\overset{\bullet}{\|}_{1/2,2,1}$ is bounded by a fixed constant depending only on $\gamma$. Therefore $K_h^{\gamma}\in V_2(\R)$ and the entropy bound (\ref{deux}) is obtained by applying  Lemma 1 of \cite{GN2010-aop}.

\subsection{\textbf{Proof of Lemma~\ref{Lepski}}}
\label{proofLepski}

\noindent We recall that the bandwidth $h_{\mathbf{\widehat m }}$ with $\mathbf{\widehat m }$ is defined in (\ref{Lep-proc2}). Let $r_n(x) = \max\left( \sqrt{\frac{1+x}{n}} , \frac{1+x}{n}\right)$ and define
\begin{eqnarray}
        \label{m*}
       m^{*}: = \mathrm{argmin}_{1\leq m \leq M} \left\{ h_m^{r/2-1} e^{-\frac{\beta}{h_m^{r}}} 
       + e^{\gamma h_{m}^{-2}}r_n(x+\log M) \right\},
\end{eqnarray}
\noindent and
\begin{eqnarray*}
 B(m) = \max_{j:j>m}\left\{  \Vert \widehat{W}^\gamma_{h_{m}} - \widehat{W}^\gamma_{h_{j}}\Vert_{\infty} - 2 \kappa e^{\gamma h_j^{-2}} r_n(x+\log M) \right\}.
\end{eqnarray*}
\noindent In one hand, we have
\begin{eqnarray*}\hspace{-2cm}
\Vert \widehat{W}^\gamma_{h_{\mathbf{\widehat m }}} - \widehat{W}^\gamma_{h_{m^*}}\Vert_{\infty} \1_{\mathbf{\widehat m } > m^*}  \hspace{-0.15cm}&=&\hspace{-0.15cm} \left(\Vert \widehat{W}^\gamma_{h_{\mathbf{\mathbf{\widehat m } }}} - \widehat{W}^\gamma_{h_{m^*}}\Vert_{\infty} -2\kappa e^{\gamma h_{\mathbf{\widehat m }}^{-2}}r_n(x+\log M)\right)\1_{\mathbf{\widehat m } > m^*}\\
&& + 2\kappa e^{\gamma h_{\mathbf{\widehat m }}^{-2}}r_n(x+\log M)\1_{\mathbf{\widehat m } > m^*}\\
\hspace{-0.15cm}&\leq&\hspace{-0.15cm} \left(B(m^*) + 2\kappa e^{\gamma h_{\mathbf{\widehat m }}^{-2}}r_n(x+\log M)\right)\1_{\mathbf{\widehat m } > m^*}.
\end{eqnarray*}
\noindent In the other hand, similarly, we have
\begin{eqnarray*}\hspace{-2cm}
\Vert \widehat{W}^\gamma_{h_{\mathbf{\widehat m }}} - \widehat{W}^\gamma_{h_{m^*}}\Vert_{\infty} \1_{\mathbf{\widehat m } \leq m^*}  &\hspace{-0.15cm}\leq&\hspace{-0.15cm}  \left(B(\mathbf{\widehat m }) + 2\kappa e^{\gamma h_{m^*}^{-2}}r_n(x+\log M)\right)\1_{\mathbf{\widehat m } \leq m^{*}}.
\end{eqnarray*}
\noindent Combining the last two displays, and by definition of $\Lcal_\kappa(\cdot)$ in (\ref{Lep-proc1}), we get
\begin{eqnarray}
       \hspace{-2cm}
\Vert \widehat{W}^\gamma_{h_{\mathbf{\widehat m }}} - \widehat{W}^\gamma_{h_{m^*}}\Vert_{\infty} \hspace{-0.15cm}&\leq&\hspace{-0.15cm} \left(B(m^*) + 2\kappa e^{\gamma h_{\mathbf{\widehat m }}^{-2}}r_n(x+\log M)\right)\1_{\mathbf{\widehat m } > m^*}\nonumber\\
&& + \left(B(\mathbf{\widehat m }) + 2\kappa e^{\gamma h_{m^*}^{-2}}r_n(x+\log M)\right)\1_{\mathbf{\widehat m } \leq m^*}\nonumber\\
\hspace{-0.15cm}&\leq&\hspace{-0.15cm}B(m^*) + B(\mathbf{\widehat m }) + 2\kappa r_n(x+\log M)(e^{\gamma h_{\mathbf{\widehat m }}^{-2}}+e^{\gamma h_{m}^{-2}})\nonumber\\
\label{a6}
\hspace{-0.15cm}&=&\hspace{-0.15cm}\mathcal L(m^*) + \mathcal L(\mathbf{\widehat m })\leq 2 \mathcal L( m^*),
\end{eqnarray}
\noindent where the last inequality follows from the definition of $\mathbf{\widehat m }$ in (\ref{Lep-proc2}).  By the definition of $B(\cdot)$, it comes 
\begin{eqnarray*}\hspace{-2cm}
\mathcal L( m^*) \hspace{-0.15cm}&=&\hspace{-0.15cm}B(m^*) + 2 \kappa e^{\gamma h_{m^*}^{-2}} r_n(x+\log M)\\
\hspace{-0.15cm}&=&\hspace{-0.15cm} \max_{j:j>m^*}\left\{  \Vert \widehat{W}^\gamma_{h_{m^*}} - \widehat{W}^\gamma_{h_{j}}\Vert_{\infty} - 2 \kappa e^{\gamma h_j^{-2}} r_n(x+\log M) \right\}\\
\hspace{-0.15cm}&\leq&\hspace{-0.15cm} \max_{j:j>m^*} \left\{ \Vert\widehat{W}^\gamma_{h_{m^*}} - \E[\widehat{W}^\gamma_{h_{m^*}}] \Vert_{\infty}+ \Vert \E[\widehat{W}^\gamma_{h_{m^*}}] -W_\rho\Vert_{\infty} + \Vert W_\rho - \E[\widehat{W}^\gamma_{h_{j}}]\Vert_{\infty} \right.\\
  &&\left.+ \Vert\E[\widehat{W}^\gamma_{h_j} ]- \widehat{W}^\gamma_{h_{j}}\Vert_{\infty}] - 2 \kappa e^{\gamma h_j^{-2}} r_n(x+\log M) \right\}+ 2 \kappa e^{\gamma h_{m^*}^{-2}} r_n(x+\log M).
  \end{eqnarray*}  %

\noindent  On the event $\mathcal E_{\kappa}$, it follows that
\begin{eqnarray*}\hspace{-2cm}
\mathcal L( m^*) \hspace{-0.15cm}&\leq&\hspace{-0.15cm}  \max_{j:j>m^*} \left\{ \Vert\widehat{W}^\gamma_{h_{m^*}} - \E[\widehat{W}^\gamma_{h_{m^*}}] \Vert_{\infty}+ \Vert \E[\widehat{W}^\gamma_{h_{m^*}}] -W_\rho\Vert_{\infty} + \Vert W_\rho - \E[\widehat{W}^\gamma_{h_{j}}]\Vert_{\infty} \right.\\
  &&\left.- \kappa e^{\gamma h_j^{-2}} r_n(x+\log M) \right\}+ 2 \kappa e^{\gamma h_{m^*}^{-2}} r_n(x+\log M).
 \end{eqnarray*}
\noindent   As $h_{m^*}>h_{j}$ for all $j>m^*$, we have $ - e^{\gamma h_{j}^{-2}} <-e^{\gamma h_{m^*}^{-2}} $. Therefore, on the event $\mathcal E_{\kappa}$, we get
\begin{eqnarray}
       \hspace{1cm}
\label{a4}
\mathcal L( m^*)  \hspace{-0.15cm}&\leq&\hspace{-0.15cm} \Vert \E[\widehat{W}^\gamma_{h_{m^*}}] -W_\rho\Vert_{\infty}+  \max_{j:j>m^*} \left\{ \Vert \E[\widehat{W}^\gamma_{h_{j}}]- W_\rho \Vert_{\infty}\right\}%
+ 2 \kappa e^{\gamma h_{m^*}^{-2}} r_n(x+\log M).
\end{eqnarray}
\noindent From (\ref{a6}) and on the event  $\mathcal E_{\kappa}$, we have
\begin{eqnarray*}\hspace{-2cm}
\label{a5}
\|\widehat{W}^\gamma_{h_{\mathbf{\widehat m }}}  - W_{\rho}\|_{\infty}\hspace{-0.15cm}&\leq&\hspace{-0.15cm} \|\widehat{W}^\gamma_{h_{\mathbf{\widehat m }}}  - \widehat{W}^\gamma_{h_{m^*}}\|_{\infty}+\|\widehat{W}^\gamma_{h_{m^*}} - W_{\rho}\|_{\infty}%
\leq|\widehat{W}^\gamma_{h_{m^*}} - W_{\rho}\|_{\infty}+2\mathcal L( m^*) \nonumber\\
\hspace{-0.15cm}&\leq&\hspace{-0.15cm}\|\widehat{W}^\gamma_{h_{m^*}} - \E[\widehat{W}^\gamma_{h_{m^*}}]\|_{\infty}+\|\E[\widehat{W}^\gamma_{h_{m^*}}] - W_{\rho}\|_{\infty}+2\mathcal L( m^*) \nonumber\\
\hspace{-0.15cm}&\leq&\hspace{-0.15cm}\kappa e^{\gamma h_{m}^{-2}}r_n(x+\log M)+\|\E[\widehat{W}^\gamma_{h_{m^*}}] - W_{\rho}\|_{\infty}+2\mathcal L( m^*). 
\end{eqnarray*}
\noindent Combining the last inequality with (\ref{a4}) 
\begin{eqnarray*}%
\|\widehat{W}^\gamma_{h_{\mathbf{\widehat m }}}  - W_{\rho}\|_{\infty}\hspace{-0.15cm}&\leq&\hspace{-0.15cm} 5\kappa e^{\gamma h_{m}^{-2}}r_n(x+\log M)+3|\E[\widehat{W}^\gamma_{h_{m^*}}] - W_{\rho}\|_{\infty}%
+2 \max_{j:j>m^*} \left\{ \Vert \E[\widehat{W}^\gamma_{h_{j}}]- W_\rho \Vert_{\infty}\right\}.
\end{eqnarray*}
\noindent From Proposition~\ref{bias}, the bias is bounded by $t\rightarrow t^{r/2-1}e^{-\beta t^{-r}}$ an increasing function for sufficiently small $t>0$, and as s $h_{m^*}>h_{j}$ for all $j>m^*$, we can write

\begin{eqnarray*}\hspace{-2cm}
\|\widehat{W}^\gamma_{h_{\mathbf{\widehat m }}}  - W_{\rho}\|_{\infty}\hspace{-0.15cm}&\leq&\hspace{-0.15cm}C\left(\kappa e^{\gamma h_{m}^{-2}}r_n(x+\log M)+h_{m^*}^{r/2-1} e^{-\beta h_{m^*}^{-r}}  \right).\nonumber\\
\end{eqnarray*}
\noindent The result comes from (\ref{m*}), the definition of $m^*$.

\subsection{\textbf{Proof of Lemma~\ref{fact1}}}
\label{prooffact1}

\noindent In view of Fatou's Lemma, we have
\begin{eqnarray*}\hspace{-1cm}
                 \mathrm{liminf}_{|z|\rightarrow \infty} z^2 p_0^{\gamma}(z) \hspace{-0.15cm}
                 &\geq&\hspace{-0.15cm} \int  \mathrm{liminf}_{|z|\rightarrow \infty} z^2 
                  p_0(z-x) N^{\gamma}(x)dx\\
                  \hspace{-0.15cm}&\geq&\hspace{-0.15cm}\int_{-\sqrt{2\gamma}}
                  ^{\sqrt{2\gamma}} \mathrm{liminf}_{|z|\rightarrow \infty} z^2 p_0(z-x)
                   N^{\gamma}(x) dx.
\end{eqnarray*}
Recall that $\gamma = \frac{1-\eta}{4\eta} \leq 1/4$, then for $|z|\geq \sqrt{2\gamma}+1$ and any $x\in (-\sqrt{2\gamma}, \sqrt{2\gamma})$, it comes by  Lemma \ref{p0ltailbound} that $p_0(z-x) \geq c (z-x)^{-2}$. Thus,
\begin{eqnarray*}\hspace{-1cm}
\mathrm{liminf}_{|z|\rightarrow \infty} z^2 p_0^{\gamma}(z)\geq c \int_{-\sqrt{2\gamma}}^{\sqrt{2\gamma}} N^{\gamma}(x) dx =  c \int_{-1}^{1} \frac{1}{\sqrt{2\pi}} e^{-\frac{x^2}{2}}dx \geq c'>0,
\end{eqnarray*}
\noindent where $c'>0$ is a numerical constant . Choose now a numerical constant $\widetilde c\geq 0$ such that $\int_{-\widetilde c}^{\widetilde c} p_{0}(x)dx\geq 1/2$, therefore, for any $|z|\leq 1+\sqrt{2\gamma}$ and some numerical constant $c''>0$ we get
\begin{eqnarray*}\hspace{-1cm}
          p_0^{\gamma}(z) \hspace{-0.15cm}
                 &\geq&\hspace{-0.15cm} \int_{-\widetilde c}^{\widetilde c} p_{0}(x)N^{\gamma}(z-x)dx \geq \min_{|y|\leq M+1 +\sqrt{2\gamma}}\{N^{\gamma}(y)\}\int_{-\widetilde c}^{\widetilde c} p_{0}(x)dx\\
                 \hspace{-0.15cm}
                 &\geq&\hspace{-0.15cm}\frac{1}{2}\min_{|y|\leq M+1 +\sqrt{2\gamma}}\{N^{\gamma}(y)\} \geq c''>0.
\end{eqnarray*}

\subsection{\textbf{Lemma~\ref{bonne matrice}}}
\label{proofbonne matrice}

\begin{Lemma}
\label{bonne matrice} 
\noindent The density matrix $\rho^{(m,h)}$ defined in (\ref{bienMD}) satisfies
the following conditions are satisfied :
\begin{enumerate}
     \item[(i)] Self adjoint: $\rho^{(m,h)}=(\rho^{(m,h)})^{*}$.
     \item[(ii)]  Positive semi-definite: $\rho^{(m,h)}\geq 0$.
     \item[(iii)]  Trace one: $\mathrm{Tr}(\rho^{(m,h)})=1$.
\end{enumerate}
\end{Lemma}

\noindent \textbf{Proof:}
\noindent \\
\noindent $\bullet$  Note first that $V_{m,h}$ is not a Wigner function, however it belongs to the linear spans of Wigner functions. Consequently, it admits the following representation
$$
\frac{1}{\pi}\mathcal R[V_{m,h}] (x,\phi) \1_{(0,\pi}(\phi)= \sum_{j,k=0}^{\infty} \tau_{j,k}^{(m,h)} \psi_j(x) \psi_{k}(x)e^{-i(j-k)\phi},
$$
\noindent where
 \begin{equation}
 \label{taujk}
\tau_{j,k}^{(m,h)} =\int_{0}^{\pi}\hspace{-0.25cm}\int\frac{1}{\pi}\mathcal R[V_{m,h}] (x,\phi)\psi_j(x) \psi_{k}(x)e^{-i(j-k)\phi} dx d\phi.
\end{equation}
\noindent For the sake of brevity, we set from now on $\tau = \tau^{(m,h)}$. Note that the matrix $\rho^{(m,h)}$ satisfies $\rho^{(m,h)}_{j,k} = \rho_{j,k}^{(0)} + \tau_{j,k}$. Exploiting the above representation of $\tau$, it is easy to see that $\tau_{j,k} = \overline{\tau_{k,j}}$ for any $j,k\geq 0$. On the other hand, $\rho^{(0)}$ is a diagonal matrix with real-valued entries. This gives (i) immediately.\\

\noindent $\bullet$ We consider now (iii). First, note that $\mathcal R[V_{m,h}](\cdot,\phi)$ is an odd function for any fixed $\phi$. Indeed, its Fourier transform with respect of the frist variable 
$$
\mathcal F_1 \left[\mathcal R[V_{m,h}](\cdot,\phi)\right] (t) = \widetilde{V}_{m,h}(t\cos \phi,t\sin \phi),
$$  is an odd function of $t$ for any fixed $\phi$. Thus, it is easy to see that $\tau_{j,j} = 0$, for any $j\geq 0$. Since $\rho^{(0)}$ is already known to be a density matrix, this implies that
$$
\mathrm{Tr}(\rho^{(m,h)}) = \mathrm{Tr}(\rho^{(0)}) + \mathrm{Tr}(\tau) = 1.
$$ 

\noindent $\bullet$ Now prove (ii). From (\ref{pattljk}), we have 
$$
|\widetilde f_{k,j}(t)| = \pi^2 |t|l_{j,k} (|t|/2),\quad j\geq k.
$$
Moreover by Lemma~\ref{ljk},  we have
 \begin{equation*}
     l_{j,k}(x) \leq \frac1\pi\left\{ \begin{array}{ll}
                                 1 &{\rm \ if \ } 0 \leq x \leq \sqrt{j+k+1}, \\
                                 e^{-(x-\sqrt{j+k+1})^2} &{\rm \ if \ } x\geq \sqrt{j+k+1}.
                                       \end{array} 
                               \right. 
     \end{equation*}
\noindent Therefore, by the change of variable ($t,\phi)$ into $w=(w_1,w_2)$, (\ref{taujk}) is such that 
\begin{eqnarray}\hspace{-1cm}
\label{I1I2}
|\tau_{j,k}| \hspace{-0.15cm}&\leq&\hspace{-0.15cm} \frac{1}{\pi} \int_{0}^{\pi}\hspace{-0.25cm}\int |\widetilde V(t\cos \phi, t \sin \phi)| |\widetilde f (t)|dt%
=\pi \iint |\widetilde V(w)| l_{j,k}\left( \frac{\|w\|}{2}  \right) dw\nonumber\\
\hspace{-0.15cm}&\leq&\hspace{-0.15cm} \int_{\|w\| \leq \sqrt{J}} |\widetilde V(w)|dw + \int_{\|w\| > \sqrt{J}} |\widetilde V(w)| e^{-(\|w\| - J)^2} dw = I_1 + I_2,
\end{eqnarray}
where $J = j+k+1$. The term  $I_1$ can be bounded as follow
\begin{eqnarray*}\hspace{-1cm}
I_1 \hspace{-0.15cm}&=&\hspace{-0.15cm} a C_0 h^{-1} e^{\beta h^{-2}} \int_{\|w\| \leq \sqrt{J}} e^{- 2  \beta \|w\|^2 } |g_{m}(\|w\|^2 - h^{-2}) g(w_2)|dw\\
\hspace{-0.15cm}&\leq&\hspace{-0.15cm}a c_0 h^{-1} e^{\beta h^{-2}} \int_{\|w\|^2 \leq J} e^{-2\beta \|w\|^2} |g_m(\|w\|^2 - h^{-2})|dw,
\end{eqnarray*}
where $C_0 = \sqrt{\pi L (\beta+\gamma)}$.

If $k+j +1 < a_m^2$, then $I_1 = 0$. If $k+j+1\geq a_m^2$, then
\begin{eqnarray}\hspace{-1cm}
\label{I1a}
I_1 \hspace{-0.15cm}&\leq&\hspace{-0.15cm} a C_0 h^{-1} e^{\beta h^{-2}} \iint  e^{-2\beta \|w\|^2}1_{a_m\leq \|w\| \leq b_m} dw\leq%
a \delta C_0 h^{-1}\delta^2 e^{\beta h^{-2}} e^{-2\beta a_m^2}\nonumber\\
\hspace{-0.15cm}&\leq&\hspace{-0.15cm}a \delta C_0 \delta^2 h^{-1}\delta^2 e^{-\beta a_m^2} \leq a C_1 \delta e^{-\beta J},
\end{eqnarray}
where $C_1>0$ is a constant depending only on $L,\beta,\gamma$.  Similarly for $I_2$, we get
\begin{eqnarray*}\hspace{-1cm}
             I_2 \hspace{-0.15cm}&=&\hspace{-0.15cm}a C_0 h^{-1} e^{\beta h^{-2}}
              \int_{\|w\| \geq \sqrt{J}} e^{- 2  \beta \|w\|^2 } |g_{m}(\|w\|^2 - h^{-2}) g(w_2)|
               e^{(\|w\| - \sqrt{J})^2} dw\\
               \hspace{-0.15cm}&\leq&\hspace{-0.15cm} a C_0 h^{-1} e^{\beta h^{-2}} 
               \int_{\|w\|^2 \geq J} e^{-2\beta \|w\|^2} |g_m(\|w\|^2 - h^{-2})| e^{(\|w\| - \sqrt{J})^2}  
               dw.
\end{eqnarray*}
If $k+j +1 \geq  b_m^2$, then $I_2 = 0$, otherwise if $k+j+1\leq b_m^2$, we have
\begin{eqnarray}\hspace{-1cm}
\label{I2a}
                   I_2 \hspace{-0.15cm}&\leq&\hspace{-0.15cm}a C_0 h^{-1} e^{\beta h^{-2}}
                    \iint  e^{-2\beta \|w\|^2}1_{a_m\leq \|w\| \leq b_m} e^{(\|w\| - \sqrt{J})^2} dw
                    \nonumber\\
                    \hspace{-0.15cm}&\leq&\hspace{-0.15cm}a \delta C_0' h^{-1}\delta^2
                     e^{\beta h^{-2}} e^{-2\beta a_m^2}\leq a \delta C_1' h^{-1} 
                     \delta^2 e^{-\beta a_m^2}\nonumber\\
                      \hspace{-0.15cm}&\leq&\hspace{-0.15cm}a \delta C_1' h^{-1} \delta^2 
                      e^{-\beta b_m^2} e^{\beta \delta}  \leq a C_1\delta e^{\beta \delta} 
                      e^{-\beta J}.
\end{eqnarray}
Combining (\ref{I1I2}), (\ref{I1a}) and (\ref{I2a}), we get for any $j\neq k$ that
\begin{eqnarray*}\hspace{-1cm}
|\tau_{j,k}| \leq c  a \delta e^{\beta \delta} e^{-\beta (j+k+1)},
\end{eqnarray*}
for some numerical constant $c>0$. Since $\rho$ is an Hermitian matrix (\textit{iii}), it admits real eigenvalues. For any eigenvalue $\lambda$ of $\rho$, in view of Theorem \ref{Gershgorin} below, there exists an integer $j\geq 1$ such that
\begin{eqnarray}
       \hspace{-1cm}\label{taubound}
\left| \lambda - \rho_{jj}^{(0)} \right| \hspace{-0.15cm}&\leq&\hspace{-0.35cm}\sum_{k =1\,:\, k\neq j}^{\infty} |\tau_{jk}| \leq c a \delta \frac{e^{\beta \delta - 2\beta}}{1-e^{-\beta}} e^{-\beta j}=: r_{j}.
\end{eqnarray}
Recall that $\rho^{(0)} = \rho^{\alpha,\lambda}$ for some $0<\alpha,\lambda<1$ where $\rho^{\alpha,\lambda}$ is defined in (\ref{rhoalphalambda}). Lemme 2 in the paper of \cite{Butucea&Guta&Artiles} guarantees that 
$$
\rho_{jj}^{\alpha,\lambda} = \frac{\alpha}{(1-\lambda)^{\alpha}} \Gamma(\alpha + 1) j^{-(1+\alpha)} (1+o(1)),
$$
as $n\rightarrow \infty.$ We note that $\rho_{jj}^{(0)}>0$ decreases polynomially with $j$ whereas $r_{j}$ decreases exponentially. Taking the numerical constant $a>0$ small enough in (\ref{jh}) independently of $j$, we get $\rho_{jj} \geq \frac{r_j}{2}\geq 0$. Thus $\rho$ is positive semi-definite.

\begin{Theorem}[Gershgorin Disk Theorem]\label{Gershgorin}
Let $A$ be an infinite square matrix and let $\mu$ be any eigenvalue of $A$. Then, for some $j\geq 1$, we have
$$
|\mu - A_{j,j}| \leq r_{j}(A),
$$ 
where $r_j(A) = \sum_{k\geq 1: k\neq j} |A_{j,k}|$.
\end{Theorem}

\noindent\textbf{{Proof: }}
Let $\mu$ be an eigenvalue of $A$ with associated unit eigenvector $v = (v_1,v_2,\ldots)$. We have
\begin{eqnarray*}\hspace{-1cm}
\lambda v_k = [A v]_k = \sum_{l\geq 1} A_{kl} v_l.
\end{eqnarray*}
We set $\widetilde k = \mathrm{argmax}_{k\geq 1}(|v_{k}|)$. Then
$$
(\mu - A_{\widetilde{k}\widetilde{k}})v_{\widetilde{k}} = \sum_{l\,:\, l\neq \widetilde{k}} A_{\widetilde{k}l} v_l.
$$
Consequently
$$
|\mu - A_{\widetilde{k}\widetilde{k}}| \leq  \sum_{l\,:\, l\neq \widetilde{k}}| A_{\widetilde{k}l} | \frac{|v_l|}{|v_{\widetilde k}|}\leq \sum_{l\,:\, l\neq \widetilde{k}}| A_{\widetilde{k}l} | := r_{\widetilde{k}}(A).
$$

\section*{Acknowledgement}
\noindent The work of K. Lounici is supported by Simons Collaboration Grant 315477 and NSF CAREER Grant DMS-1454515. The work of K. Meziani is supported by  "Calibration" ANR-2011-BS01-010-01. The work of G. Peyr\'e is supported by the European Research Council (ERC project SIGMA-Vision). 
\bibliographystyle{imsart-nameyear}
\bibliography{bibli}

\begin{thebibliography}{100}
{\small
\providecommand{\natexlab}[1]{#1}
\providecommand{\url}[1]{\texttt{#1}}
\expandafter\ifx\csname urlstyle\endcsname\relax
  \providecommand{\doi}[1]{doi: #1}\else
  \providecommand{\doi}{doi: \begingroup \urlstyle{rm}\Url}\fi
  
  
  
  
\bibitem[Alquier, Meziani and Peyr\'e(2013)]{AMP}
  \textsc{Alquier, P., Meziani, K. and Peyr\'e,G.},
  \newblock \emph{ Adaptive Estimation of the Density Matrix in Quantum Homodyne Tomography with Noisy Data}.
\newblock {Inverse Problems}, 29, 7, 075017, 2013.





 \bibitem[Aubry, Butucea and Meziani(2009)]{ABM}
\textsc{Aubry, J.-M. and Butucea, C. and Meziani, K.},
 \newblock \emph{State estimation in quantum homodyne tomography with noisy data}.
\newblock{Inverse Problems}, 25, 1, 2009.
  
  
  
  
  
  
 \bibitem[Artiles, Gill and Gu\c{t}\u{a}(2005)] {Artiles&Gill&Guta}
\textsc{Artiles, L. and Gill, R. and Gu{\c{t}}{\u{a}}, M.},
\newblock\emph{An invitation to quantum tomography}.
\newblock{J. Royal Statist. Soc. B (Methodological)}, 67,109--134, 2005.







\bibitem[{A}verbuch et al.(2008)]{SlantStack}
\textsc{{A}verbuch, A., {C}oifman, R.R., {D}onoho, D.L., {I}sraeli, {M}., {S}hkolnisky, {Y}. and {S}edelnikov, {I}.} ,
\newblock \emph{{A} framework for discrete integral transformations: {I}{I}. {T}he 2{D} discrete {R}adon transform.}.
\newblock \emph{{S}{I}{A}{M} {J}. {S}ci. {C}omput.}, 30(2), 785--803, 2008.



\bibitem[Barndorff-Nielsen, Gill and Jupp(2003)]{Barndorff-Nielsen&Gill&Jupp}
\textsc{Barndorff-Nielsen, O.~E. and Gill, R. and Jupp, P.~E.},
\newblock \emph{On quantum statistical inference (with discussion)}.
\newblock {J. Royal Stat. Soc. B}, 65, 775--816, 2003.



\bibitem[Bergh and L{\"o}fstr{\"o}m(1976)]{BerghLof76}
\textsc{Bergh, J. and L{\"o}fstr{\"o}m, J.},
\newblock \emph{Interpolation spaces. {A}n introduction}.
\newblock {Grundlehren der Mathematischen Wissenschaften, No. 223,Springer-Verlag, Berlin}, 1976.


\bibitem[Bousquet(2002)]{Bousquet2002}
    \textsc{Bousquet, O.},
    \newblock\emph{A {B}ennett concentration inequality and its application to
              suprema of empirical processes}.
 \newblock{C. R. Math. Acad. Sci. Paris}, 334, 6, 495--500, 2002.




\bibitem[Bourdaud, Lanza de Cristoforis and Sickel(2006)] {Bourdaudal2006}
\textsc {Bourdaud, G. and Lanza de Cristoforis, M. and Sickel, W.},
\newblock \emph{Superposition operators and functions of bounded
              {$p$}-variation}.
\newblock \emph{Rev. Math. Iberoam.}, 2, 455--487, 2006.




\bibitem[Butucea, Gu{\c{t}}{\u{a}} and Artiles(2007)]{Butucea&Guta&Artiles}
\textsc{Butucea, C. and Gu{\c{t}}{\u{a}}, M. and Artiles, L.},
\newblock\emph{Minimax and adaptive estimation of the Wigner
                  function in quantum homodyne tomography with noisy
                  data}.
\newblock{Ann. Statist.}, 2, 35, 465--494, 2007.   

   \bibitem[D'Ariano, Macchiavello and Paris(1994)]{DAriano.0}
\textsc{D'Ariano, G.~M. and Macchiavello, C. and Paris, M.~G.~A.},
\newblock \emph{Detection of the density matrix through optical
                  homodyne tomography without filtered back
                  projection}.
\newblock {Phys. Rev. A}, 50, 4298--4302, 1994.





\bibitem[Erd\'elyi , Magnus, Oberhettinger and Tricomi(1953)]{Erdelyi}
\textsc{Erd\'elyi A., Magnus W., Oberhettinger F., Tricomi F.G.}, 
\newblock \emph{Higher transcendental functions}.
\newblock {McGraw-Hill Book Company, Inc., New York-Toronto-
London }, Vols. I, II, 1953. 





\bibitem[Gin{\'e} and Nickl(2009)]{GN2010-aop}
   \textsc{Gin{\'e}, E. and Nickl, R.},
    \newblock\emph{Uniform limit Theorems for wavelet density estimators}.
  \newblock{Ann. Probab.}, 37, 4,1605--1646, 2009.
 
 
 
 
  \bibitem[Gu{\c{t}}{\u{a}} and Artiles(2007)]{Guta&Artiles}
\textsc{Gu{\c{t}}{\u{a}}, M. and Artiles, L.},
\newblock\emph{Minimax estimation of the Wigner
           in quantum homodyne tomography with ideal detectors}.
\newblock{Math. Methods Statist.}, 16, 1,1--15, 2007.



   
  
  
\bibitem[Helstrom(1976)]{Helstrom}
\textsc{Helstrom, C.~W.},
\newblock \emph{Quantum Detection and Estimation Theory}.
\newblock {Academic Press, New York}, 1976.



\bibitem[Holevo(1982)]{Holevo}
\textsc{Holevo, A.~S.},
\newblock \emph{Probabilistic and Statistical Aspects of Quantum Theory}.
\newblock {North-Holland}, 1982.

\bibitem[Leonhardt(1997)]{Leonhardt}
\textsc{Leonhardt, U.},
\newblock \emph{Measuring the Quantum State of Light}.
\newblock {Cambridge University Press}, 1997.



\bibitem[Leonhardt, Paul and D'Ariano(1995)]{DAriano.3}
\textsc{Leonhardt, U. and Paul, H. and D'Ariano, G.~M.},
\newblock \emph{Tomographic reconstruction of the density matrix via
                  pattern functions}.
\newblock {Phys. Rev. A}, 52, 4899--4907, 1995.


 \bibitem[Lounici, K. and Nickl(2001)]{LN2011}
  \textsc{Lounici, K. and Nickl, R.},
  \newblock \emph{Global uniform risk bounds for wavelet deconvolution
              estimators}.
\newblock {Ann. Statist.}, 39, 1, 201--231, 2001.


\bibitem[Meziani(2008)]{MezianiTest}
 \textsc{Meziani, K.},
\newblock \emph{Nonparametric goodness-of fit testing in quantum homodyne
              tomography with noisy data}.
\newblock{Electron. J. Stat.}, 2,  1195--1223, 2008.




\bibitem[Meziani(2007)]{Meziani}
 \textsc{Meziani, K.}
\newblock \emph{Nonparametric Estimation of the Purity of a Quantum
                  State in Quantum Homodyne Tomography with Noisy
                  Data}.
\newblock{Math. Meth. of Stat.}, 4, 16, 1--15, 2007.

\bibitem[Richter(2000)]{Richter1}
\textsc{Richter, T.},
\newblock \emph{Realistic pattern functions for optical homodyne
                  tomography and determination of specific expectation
                  values}.
\newblock {Phys. Rev. A}, 61, 2000.



\bibitem[Vogel and Risken(1989)]{Vogel&Risken}
\textsc{Vogel, K. and Risken, H. },
\newblock \emph{Determination of quasiprobability distributions in
                  terms of probability distributions for the rotated
                  quadrature phase}.
\newblock {Phys. Rev. A}, 40, 2847--2849, 1989.


\bibitem[Wigner(1932)]{Wigner}
\textsc{Wigner, E.},
\newblock \emph{On the quantum correction for thermodynamic
                  equations}.
\newblock {Phys. Rev.}, 40, 749--759, 1932.




\bibitem[Tsybakov(2009)]{Tsybakovlivre}
\textsc{Tsybakov, A.B.},
\newblock \emph{Introduction to Nonparametric Estimation}.
\newblock {Springer Series in Statistics, New York}, 2009.


}
\end{thebibliography}

\end{document}